%% file: LKS-CUT0.tex
\documentclass[A4paper,oneside,11pt]{article}
\usepackage{amssymb,amsthm,amsmath,marvosym,amsfonts,fontenc}
\usepackage[usenames,dvipsnames]{color}
\usepackage{enumerate,marginnote}
\usepackage[dvips]{graphicx}
\usepackage{mdwlist}
\usepackage{pifont}
\usepackage{bm}
\usepackage{multind}\makeindex{general}\makeindex{mathsymbols}
\usepackage{caption}
\usepackage{fullpage} 
\usepackage{epsfig,psfrag}
\usepackage{fancyhdr} 
\usepackage{nohyperref}

\usepackage{xr}
\numberwithin{equation}{section}
\numberwithin{table}{section} 
\numberwithin{figure}{section}
\newtheorem{theorem}{Theorem}[section]

\newtheorem{lemma}[theorem]{Lemma}
\newtheorem{conjecture}[theorem]{Conjecture}

\newtheorem{proposition}[theorem]{Lemma}
\newtheorem{cor}[theorem]{Corollary}

\newtheorem{fact}[theorem]{Fact}

\newtheorem{definition}[theorem]{Definition}

\theoremstyle{remark}

\newcommand{\By}[2]{\overset{\mbox{\tiny{#1}}}{#2}}
\newcommand{\ByRef}[2]{   \By{\eqref{#1}}{#2} }

\newcommand{\leBy}[1]{    \By{#1}{\le} }

\newcommand{\leByRef}[1]{ \ByRef{#1}{\le} }

\let\sm\setminus
\let\subset\subseteq 
 
\let\epsilon\varepsilon
\let\sharp\#
\def\dcup{\dot\cup} 
\renewcommand{\leq}{\leqslant}
\renewcommand{\le}{\leqslant}
\renewcommand{\geq}{\geqslant}
\renewcommand{\ge}{\geqslant}

\let\oldmarginpar\marginpar
\renewcommand\marginpar[1]{\-\oldmarginpar[\raggedleft\footnotesize #1]%
{\raggedright\footnotesize #1}}

\newcommand{\HIDDENPROOF}[1]{
}

\newcommand{\HIDDENTEXT}[1]{
}

\title{The approximate Loebl--Koml\'os--S\'os Conjecture~I:\\ The sparse decomposition} 
\input{affiliations.tex}
\input{definice.tex}

\renewcommand{\today}{}
\date{}

\begin{document}
\pagenumbering{roman}
\maketitle
\begin{abstract}
In a series of four papers
we prove the following relaxation of the
Loebl--Koml\'os--S\'os Conjecture: For every~$\alpha>0$
there exists a number~$k_0$ such that for every~$k>k_0$
 every $n$-vertex graph~$G$ with at least~$(\frac12+\alpha)n$ vertices
of degree at least~$(1+\alpha)k$ contains each tree $T$ of order~$k$ as a
subgraph.

The method to prove our result follows a strategy similar to
approaches that employ the Szemer\'edi regularity lemma: we
decompose the graph~$G$, find a suitable combinatorial structure inside the decomposition, and then embed the tree~$T$ into~$G$ using this structure. Since for sparse graphs~$G$, the decomposition given by the regularity lemma is not helpful, we use a more general decomposition technique. We show that each graph can be decomposed into vertices of huge degree, regular pairs (in the sense of the regularity lemma), and two other objects each exhibiting certain expansion properties. 
In this paper, we introduce this novel decomposition technique.  In the three follow-up papers, we find a combinatorial structure suitable inside the decomposition, which we then use for embedding the tree.
\end{abstract}

\bigskip\noindent
{\bf Mathematics Subject Classification: } 05C35 (primary), 05C05 (secondary).\\
{\bf Keywords: }extremal graph theory; Loebl--Koml\'os--S\'os Conjecture; tree embedding; regularity lemma; sparse graph; graph decomposition.

\newpage

\rhead{\today}

\newpage

\tableofcontents
\newpage
\pagenumbering{arabic}
\setcounter{page}{1}

\section{Introduction}\label{sec:intro}
\input{intro0.tex}

\section{Notation and preliminaries}\label{sec:preliminaries}

\input{prelim0.tex}

\section{Decomposing sparse graphs}\label{sec:class}
\input{Class0.tex}

%
%
%

%

\section{Acknowledgements}\label{sec:ACKN}
The work on this project lasted from the beginning of 2008 until 2014
and we are very grateful to the following institutions and funding bodies for
their support. 

\smallskip

During the work on this paper Hladk\'y was also affiliated with Zentrum
Mathematik, TU Munich and Department of Computer Science, University of Warwick. Hladk\'y was funded by a BAYHOST fellowship, a DAAD fellowship, 
Charles University grant GAUK~202-10/258009, EPSRC award EP/D063191/1, and by an EPSRC Postdoctoral Fellowship during the work on the project. 

Koml\'os and Szemer\'edi acknowledge the support of NSF grant
DMS-0902241.

Piguet has been also affiliated with the Institute of Theoretical Computer Science, Charles University in Prague, Zentrum
Mathematik, TU Munich, the Department of Computer Science and DIMAP,
University of Warwick, and the school of mathematics, University of Birmingham. Piguet acknowledges the support of the Marie Curie fellowship FIST,
DFG grant TA 309/2-1, a DAAD fellowship,
Czech Ministry of
Education project 1M0545,  EPSRC award EP/D063191/1,
and  the support of the EPSRC
Additional Sponsorship, with a grant reference of EP/J501414/1 which facilitated her to
travel with her young child and so she could continue to collaborate closely
with her coauthors on this project. This grant was also used to host Stein in
Birmingham.  Piguet was supported by the European Regional Development Fund (ERDF), project ``NTIS - New Technologies for Information Society'', European Centre of Excellence, CZ.1.05/1.1.00/02.0090.

Stein was affiliated with the Institute of Mathematics and Statistics, University of S\~ao Paulo, the Centre for Mathematical Modeling, University of Chile and the Department of Mathematical Engineering, University of Chile. She was
supported by a FAPESP fellowship, and by FAPESP travel grant  PQ-EX 2008/50338-0, also
CMM-Basal,  FONDECYT grants 11090141 and  1140766. She also received funding by EPSRC Additional Sponsorship EP/J501414/1.

We enjoyed the hospitality of the School of Mathematics of University of Birmingham, Center for Mathematical Modeling, University of Chile, Alfr\'ed R\'enyi Institute of Mathematics of the Hungarian Academy of Sciences and Charles University, Prague, during our long term visits.

The yet unpublished work of Ajtai, Koml\'os, Simonovits, and Szemer\'edi on the Erd\H{o}s--S\'os Conjecture was the starting point for our project, and our solution crucially relies on the methods developed for the Erd\H{o}s-S\'os Conjecture. Hladk\'y, Piguet, and Stein are very grateful to the former group for explaining them those techniques.

Hladk\'y would like to thank Maxim Sviridenko for discussion on the algorithmic aspects of the problem.

\medskip
A doctoral thesis entitled \emph{Structural graph theory} submitted by Hladk\'y in September 2012 under the supervision of Daniel Kr\'al at~Charles University in~Prague is based on the series of the papers~\cite{cite:LKS-cut0,cite:LKS-cut1, cite:LKS-cut2, cite:LKS-cut3}. The texts of the two works overlap greatly. We are grateful to PhD committee members Peter Keevash and Michael Krivelevich. Their valuable comments are reflected in the series. 

\bigskip

We thank the referees for their very detailed remarks.

\bigskip
The contents of this publication reflects only the authors' views and not necessarily the views of the European Commission of the European Union.

\printindex{mathsymbols}{Symbol index}
\printindex{general}{General index}

\newpage
\addcontentsline{toc}{section}{Bibliography}
\bibliographystyle{alpha}
\bibliography{bibl}
\end{document}

%% file: affiliations.tex
\author{Jan
Hladk\'y
\thanks{\emph{Corresponding author.} Institute of Mathematics, Academy of Science of the Czech Republic. \v Zitn\'a 25, 110 00, Praha, Czech Republic. The Institute of Mathematics of the Academy of Sciences of the Czech Republic is supported by RVO:67985840. Email:
\texttt{honzahladky@gmail.com}.
The research leading to these results has received funding from the People Programme (Marie Curie Actions) of the European Union's Seventh Framework Programme (FP7/2007-2013) under REA grant agreement umber 628974. Much of the work was done while supported by an EPSRC postdoctoral fellowship while affiliated with DIMAP and Mathematics Institute, University of
Warwick.}
\quad 
J\'anos Koml\'os\thanks{Department of Mathematics, Rutgers University, 110 Frelinghuysen Rd., Piscataway, NJ~08854-8019, USA} 
\quad 
Diana Piguet\thanks{Institute of Computer Science, Czech Academy of Sciences, Pod Vod\'arenskou v\v e\v z\'i 2, 182~07 Prague, Czech Republic. With institutional support RVO:67985807. Supported by the Marie Curie fellowship FIST, DFG grant TA 309/2-1,  Czech Ministry of Education project 1M0545, EPSRC award EP/D063191/1, and EPSRC Additional Sponsorship EP/J501414/1.
	The research leading to these results has received funding from the European Union Seventh
	Framework Programme (FP7/2007-2013) under grant agreement no. PIEF-GA-2009-253925.
    The work leading to this invention was supported by the European Regional Development Fund (ERDF), project ``NTIS -- New Technologies for Information Society'', European Centre of Excellence, CZ.1.05/1.1.00/02.0090.}
    \\ 
    Mikl\'os Simonovits\thanks{R\'enyi
    Institute, Budapest, Hungary. Supported by OTKA~78439, OTKA~101536, ERC-AdG.~321104} 
\quad 
Maya Stein\thanks{Department of Mathematical Engineering,
University of Chile, Santiago, Chile.  Supported by Fondecyt Iniciacion grant 11090141, Fondecyt Regular grant 1140766 and CMM Basal.}
\quad 
Endre Szemer\'edi\thanks{R\'enyi
	Institute, Budapest, Hungary. Supported by OTKA~104483 and ERC-AdG.~321104}}

%% file: definice.tex
\def\semiregular{regularized }

\newcommand{\PARAMETERPASSING}[2]{{\mathrm{#1}\ref{#2}}}

\def\NN{\mathbb{N}}

\newcommand{\V}{\mathcal V}

\def\bigboxplus{\mbox{\Large $\boxplus$}}

\def\mindeg{\mathrm{mindeg}}
\def\maxdeg{\mathrm{maxdeg}}
\def\density{\mathrm{d}}
\def\neighbour{\mathrm{N}}
\def\neighbour{\mathrm{N}}

\newcommand{\treeclass}[1]{\mathbf{trees}({#1})}

\newcommand{\LKSgraphs}[3]{\mathbf{LKS}({#1},{#2},{#3})}
\newcommand{\LKSmingraphs}[3]{\mathbf{LKSmin}({#1},{#2},{#3})}
\newcommand{\LKSsmallgraphs}[3]{\mathbf{LKSsmall}({#1},{#2},{#3})}

\newcommand{\smallvertices}[3]{\mathbb{S}_{{#1},{#2}}({#3})}
\newcommand{\largevertices}[3]{\mathbb{L}_{{#1},{#2}}({#3})}

\def\Gcapt{G_\nabla}
\def\GD{G_{\mathcal{D}}}
\def\Gblack{G_{\mathrm{reg}}}
\def\Gexp{G_{\mathrm{exp}}}
\def\BGblack{\mathbf{G}_{\mathrm{reg}}}
\def\smallatoms{\mathbb{E}}
\def\clusters{\mathbf{V}}
\def\class{\nabla}
\def\HugeVertices{\mathbb{H}}
\def\DenseSpots{\mathcal{D}}

\def\XA{\mathbb{XA}}

\def\clustersize{\mathfrak{c}}



%% file: intro0.tex
\subsection{Statement of the problem}\label{ssec:veryIntro}
This is the first of a series of four papers~\cite{cite:LKS-cut0, cite:LKS-cut1, cite:LKS-cut2, cite:LKS-cut3} in which we provide an approximate solution of the Loebl--Koml\'os--S\'os Conjecture, a problem in extremal graph theory which fits the classical form \emph{Does a certain density condition imposed on a graph guarantee a certain subgraph?} Classical results of this type include Dirac's Theorem which determines the minimum degree threshold for containment of a Hamilton cycle, or Mantel's Theorem which determines the average degree threshold for containment of a triangle. Indeed, most of these extremal problems are formulated in terms of the minimum or average degree of the host graph.

We investigate a density condition which guarantees the containment of
\emph{each} tree of order~$k$. The greedy tree-embedding strategy shows that requiring a minimum degree of more than $k-2$ is sufficient. Further, this bound
is best possible as any $(k-2)$-regular graph avoids the $k$-vertex star. Erd\H os and S\'os conjectured that one can replace the minimum degree with the average degree, with the same conclusion.
\begin{conjecture}[Erd\H os--S\'os Conjecture 1963]\label{conj:ES}
Let $G$ be a graph of average degree greater than  $k-2$. Then $G$ contains each
tree of order $k$ as a subgraph.
\end{conjecture}
A solution of the Erd\H os--S\'os Conjecture for all $k$ greater than some absolute constant was announced by Ajtai, Koml\'os, Simonovits, and Szemer\'edi in the early 1990's.
In a similar spirit, Loebl, Koml\'os, and S\'os conjectured that a \emph{median degree} of  $k-1$ or more is sufficient for containment of any tree of order $k$. By median degree we mean the degree of a vertex in the middle of the ordered degree sequence. 
\begin{conjecture}[Loebl--Koml\'os--S\'os Conjecture 1995~\cite{EFLS95}]\label{conj:LKS}
Suppose that $G$ is an $n$-vertex graph with at least $n/2$ vertices of degree more than $k-2$. Then $G$ contains each tree of order $k$.
\end{conjecture}
We discuss Conjectures~\ref{conj:ES} and~\ref{conj:LKS} in detail in
Section~\ref{ssec:LKS}. Here, we just state the main result we achieve in our series of four papers, an
approximate solution of the Loebl--Koml\'os--S\'os Conjecture.
\begin{theorem}[Main result~\cite{cite:LKS-cut3}]\label{thm:main}
For every $\alpha>0$ there exists  $k_0$ such that for any
$k>k_0$ we have the following. Each $n$-vertex graph $G$ with at least
$(\frac12+\alpha)n$ vertices of degree at least $(1+\alpha)k$ contains each tree $T$ of
order $k$.
\end{theorem}

The proof of this theorem is in~\cite{cite:LKS-cut3}. The first step towards this result is Lemma~\ref{lem:LKSsparseClass}, which constitutes the main result of the present paper. It gives a decomposition of the host graph $G$ into several parts which will be later useful for the embedding. See Section~\ref{ssec:OverviewOfProof} for a description of the result and its role in the proof of Theorem~\ref{thm:main}. 
Also see~\cite{LKS:overview} for a more detailed overview of the proof.

\subsection{The regularity lemma and the sparse decomposition}\label{ssec:iRL}
The Szemer\'edi regularity lemma has been a major tool in extremal graph theory
for more than three decades. It provides an approximation of an arbitrary graph by a collection of generalized quasi-random graphs. This allows to represent the graph by
a so-called \emph{cluster graph}.  Then, instead of solving the original problem, one can solve a modified simpler problem in the cluster graph.

The applicability of the original Szemer\'edi regularity lemma is, however, limited to \emph{dense graphs}, i.e., graphs that contain a substantial proportion of all possible edges. There is a version of the regularity lemma for sparser graphs by Kohayakawa and R\"odl~\cite{Kohayakawa97Szemeredi} later strengthened by Scott~\cite{Scott:RL}, as well as other statements that draw on something from its philosophy (e.g.~\cite{elek-2008}). However, these statements provide a picture much less informative than Szemer\'edi's original result. A regularity type representation of general (possibly sparse) graphs is one of the most important goals of contemporary discrete mathematics. By such a representation we mean an approximation of the input graph by a structure of bounded complexity carrying enough of the important information about the graph.

A central tool in the proof of Theorem~\ref{thm:main} is a structural
decomposition of the graph $G$. This
decomposition --- which we call \emph{sparse decomposition} --- applies to any
graph whose average degree is greater than a constant. The sparse decomposition provides a partition of any graph into vertices of huge degrees and into a bounded degree part. The bounded degree part is further decomposed into dense regular pairs, an edge set with certain expander-like properties, and a vertex set which is expanding in a different way (we shall give a more precise description in Section~\ref{ssec:OverviewOfProof}). 
This kind of decomposition was first used by Ajtai, Koml\'os, Simonovits, and Szemer\'edi in their yet unpublished work on the Erd\H os--S\'os Conjecture. The main goal of this paper is to present the sparse decomposition, and to show that each graph has such a sparse decomposition: This will be done in Lemma~\ref{lem:decompositionIntoBlackandExpanding}.  Lemma~\ref{lem:LKSsparseClass} provides a sparse decomposition with additional tailor-made features for graphs that fulfil the conditions of Theorem~\ref{thm:main}.

In the case of dense graphs the sparse decomposition produces a Szemer\'edi regularity partition (as explained in Section~\ref{ssec:sparsedecompofdensegraphs}), and thus the decomposition lemma (Lemma~\ref{lem:decompositionIntoBlackandExpanding}) extends the Szemer\'edi regularity lemma. 
But the interesting setting for
 the  decomposition lemma is the field of sparse graphs. 

\subsection{Loebl--Koml\'os--S\'os Conjecture and Erd\H os--S\'os Conjecture}\label{ssec:LKS}
Let us first introduce some notation. We say that $H$ \emph{embeds} in a graph $G$ and write $H\subset G$ if $H$ is a (not necessarily induced) subgraph of $G$. The associated map $\phi:V(H)\rightarrow V(G)$ is called an \index{general}{embedding}\emph{embedding of $H$ in $G$}.
More generally, for a graph class $\mathcal H$ we
write $\mathcal H\subset G$ if $H\subset G$ for
every $H\in\mathcal H$. 
Let \index{mathsymbols}{*Trees@$\treeclass{k}$}$\treeclass{k}$ be the
class of all trees of order $k$.

\medskip

Conjecture~\ref{conj:LKS} is dominated by two parameters: one quantifies the number of vertices of `large' degree, and the other tells us how large this degree should actually be. Strengthening either of these bounds sufficiently,  the conjecture becomes trivial. 
Indeed, if we replace $n/2$ with $n$, then any tree of order $k$ can be embedded greedily.
Also, if we replace $k-2$ with $4k-4$, then $G$, being a graph of average degree at least $2k-2$, has a subgraph $G'$ of minimum degree at least~$k-1$. Again we can greedily embed any tree of order $k$.

On the other hand, one may ask whether smaller lower bounds would suffice.   For the bound~$k-2$, this is not the case, since stars of order $k$ require a vertex of degree at least $k-1$ in the host graph. Another example can be obtained by considering a disjoint union of cliques of order~$k-1$. No tree of order~$k$ is contained in such a graph.  

\begin{figure}[t]
\centering 
\includegraphics[scale=0.7]{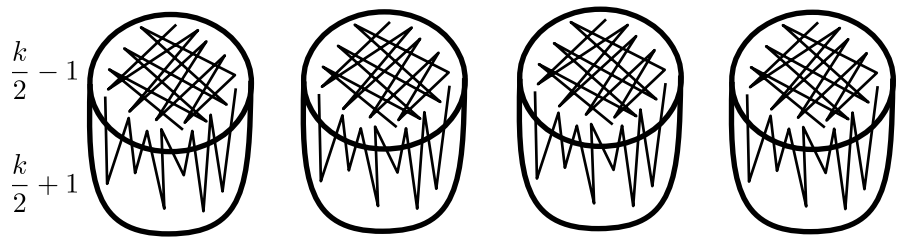}
\caption[Extremal graph for the Loebl--Koml\'os--S\'os Conjecture]{An extremal graph for the Loebl--Koml\'os--S\'os Conjecture.}
\label{fig:ExtremalGraph}
\end{figure}

For the bound $n/2$, the following example shows that this number cannot be decreased much.
First, assume that $n$ is even, and that $n=k$. Let $G^*$ be obtained from the complete graph on~$n$ vertices by deleting all edges inside a set of  $\frac
n2+1$ vertices. It is easy to
check that $G^*$ does not contain the $k$-vertex path. In general, $G^*$ does not contain any
tree of order $k$ with independence number less than $\frac k2+1$.
	 Now, taking the union of 
several disjoint copies  of $G^*$ we obtain examples for other values of $n$. (And adding a small complete component we can get to \emph{any} value of $n$.) See
Figure~\ref{fig:ExtremalGraph} for an illustration.

However, we do not know of any example attaining the exact bound
 $ n/2$. Thus it might be possible to lower the bound $n/2$ from Conjecture~\ref{conj:LKS} to the one attained in our example above:
\begin{conjecture}\label{conj:LKSstronger}
Let $k\in \NN$ and let $G$ be a graph on $n$ vertices, with more than $\frac n2-\lfloor \frac{n}{k}\rfloor - (n\mod k)$ vertices of degree at least $k-1$. Then $\treeclass{k}\subset G$.
\end{conjecture}
It might even be that
if $n/k$ is far from integrality, a slightly weaker lower bound on the number of vertices of large degree still works (see~\cite{HladkyMSC,HlaPig:LKSdenseExact}).

\medskip

Several partial results concerning 
Conjecture~\ref{conj:LKS} have been obtained; let us briefly summarize  the major ones. Two main directions can be distinguished among those results that prove the conjecture for special classes of graphs: either one places restrictions on the host graph, or on the class of trees to be embedded. Of the latter type is the result by
Bazgan, Li, and Wo{\'z}niak~\cite{BLW00}, who  proved the
conjecture for paths. Also, Piguet and Stein~\cite{PS2} proved that 
Conjecture~\ref{conj:LKS} is true for trees of
diameter at most 5, which improved earlier results of Barr and Johansson~\cite{Barr} and Sun~\cite{Sun07}.
Restrictions on the host graph have led to the following results.
Soffer~\cite{Sof00} showed that Conjecture~\ref{conj:LKS} is true if the
host graph has girth at least 7. Dobson~\cite{Dob02} proved the
conjecture for host graphs whose complement does not contain a
$K_{2,3}$. This has been extended by Matsumoto and Sakamoto~\cite{MaSa} who replace the $K_{2,3}$ with a slightly larger graph.

\medskip

A different approach is to solve the conjecture for special values of $k$. One such case, known as the Loebl conjecture, or also as the ($n/2$--$n/2$--$n/2$)-Conjecture,  is the case $k=n/2$. Ajtai, Koml\'os, and Szemer\'edi~\cite{AKS95} solved an approximate version of this conjecture, and later
 Zhao~\cite{Z07+} used a refinement of this approach to prove the sharp version of the conjecture for large
graphs.

\medskip

An approximate version of Conjecture~\ref{conj:LKS} for dense graphs, that is, for $k$ linear in $n$, was proved by Piguet
and Stein~\cite{PS07+}. 

\begin{theorem}[Piguet--Stein~\cite{PS07+}]\label{thm:PiguetStein}
For any $q>0$ and $\alpha>0$ there exists a number $n_0$ such that for any $n>n_0$ and
$k>qn$ the following holds. 
For each $n$-vertex graph $G$ with at least
$n/2$ vertices of degree at least $(1+\alpha)k$ we have $\treeclass{k+1}\subset G$.
\end{theorem}

This result was proved using the regularity method. Adding stability arguments, 
 Hladk\'y and Piguet~\cite{HlaPig:LKSdenseExact}, and independently Cooley~\cite{Cooley08} proved Conjecture~\ref{conj:LKS} for large dense graphs.
 
\begin{theorem}[Hladk\'y--Piguet~\cite{HlaPig:LKSdenseExact}, Cooley~\cite{Cooley08}]\label{thm:denseLKS}
For any $q>0$ there exists a number $n_0=n_0(q)$ such that for any $n>n_0$ and $k>qn$ the following
holds. 
For each $n$-vertex graph $G$ with at least
$n/2$ vertices of degree at least $k$ we have $\treeclass{k+1}\subset G$.
\end{theorem}

\medskip

\begin{figure}[t]
\centering 
\includegraphics[scale=0.7]{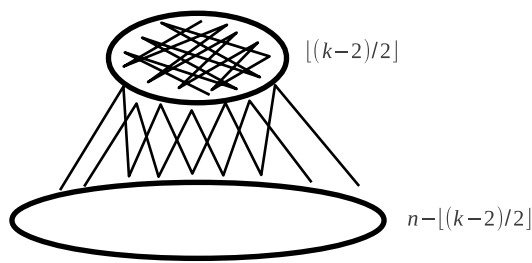}
\caption[Almost extremal graph for the Erd\H os--S\'os Conjecture]{An almost extremal graph for the Erd\H os--S\'os Conjecture.}
\label{fig:ExtremalGraphES2}
\end{figure}

Let us now turn our attention to the  Erd\H os--S\'os Conjecture. 
The Erd\H os--S\'os Conjecture~\ref{conj:ES} is best possible whenever $n(k-2)$ is even. Indeed, in that case it suffices to consider a $(k-2)$-regular graph. This is a graph with average degree exactly $k-2$ which does not contain the star of order $k$. Even when the star (which in a sense is a pathological tree) is excluded from the considerations, we can --- at least when $k-1$ divides $n$ --- consider a disjoint union of $\frac{n}{k-1}$ cliques $K_{k-1}$. This graph contains \emph{no} tree from $\treeclass{k}$.
There is another important graph with many edges which does not contain for example the path $P_{k}$, depicted in Figure~\ref{fig:ExtremalGraphES2}. 
This graph consists of a set of vertices of size $\lfloor (k-2)/2\rfloor$ that are connected to all vertices in the graph.
This graph has $\frac12(k-2)n-O(k^2)$ edges when $k$ is even and $\frac12(k-3)n-O(k^2)$ edges otherwise, and therefore gets close to the conjectured bound when $k\ll n$.
Apart from the already mentioned announced breakthrough by Ajtai, Koml\'os, Simonovits, and Szemer\'edi, work on this conjecture includes~\cite{bradob,Haxell:TreeEmbeddings,MaSa,sacwoz,woz}.

\bigskip

Both Conjectures~\ref{conj:LKS} and Conjecture~\ref{conj:ES} have an important application in Ramsey theory. Each of them implies that the Ramsey number of two trees $T_{k+1}\in\treeclass{k+1}$, $T_{\ell+1}\in\treeclass{\ell+1}$ is bounded by $R(T_{k+1}, T_{\ell+1})\le k+\ell+1$. Actually more is implied: Any $2$-edge-colouring of $K_{k+\ell +1}$ contains either \emph{all} trees in $\treeclass{k+1}$ in red, or \emph{all} trees in $\treeclass{\ell+1}$ in blue.

The bound $R(T_{k+1}, T_{\ell+1})\le k+\ell+1$ is almost tight only for certain types of trees. For example,   Gerencs\'er and Gy\'arf\'as~\cite{GerencserGyarfas} showed $R(P_k,P_\ell)=\max\{k,\ell\}+\left\lfloor\frac{\min\{k,\ell\}}2\right\rfloor-1$ for paths
$P_k\in\treeclass{k}$, $P_\ell\in\treeclass{\ell}$.
 Harary~\cite{Harary:RecentRamsey} showed $R(S_k,S_\ell)=k+\ell-2-\epsilon$ for stars $S_k\in\treeclass{k}$, $S_\ell\in\treeclass{\ell}$, where $\epsilon\in\{0,1\}$ depends on the
parity of $k$ and $\ell$. 
 Haxell, \L uczak, and Tingley confirmed asymptotically~\cite{HLT02} that the discrepancy of the Ramsey bounds for trees depends on their balancedness, at least when the maximum degrees of the trees considered are moderately bounded.

\subsection{Related tree containment problems}

\paragraph{Minimum degree conditions for spanning trees.} Recall that the tight min-degree condition for containment of a general spanning tree $T$ in an $n$-vertex graph $G$ is the trivial one, $\mindeg(G)\ge n-1$. However, the only tree which requires this bound is the star. This indicates that this threshold can be lowered substantially if we have a control of $\maxdeg(T)$. Szemer\'edi and his collaborators~\cite{KSS:SpanningTrees,CsLeNaSz:BoundedDegree} showed that this is indeed the case, and obtained tight min-degree bounds for certain ranges of $\maxdeg(T)$. For example, if $\maxdeg(T)\le n^{o(1)}$, then $\mindeg(G)\ge (\frac12+o(1))n$ is a sufficient condition. (Note that $G$ may become disconnected close to this bound.)

\paragraph{Trees in random graphs.}
To complete the picture of research involving tree containment problems we mention two rich and vivid (and also closely connected) areas: trees in random graphs, and trees in expanding graphs. The former area is centered around the following question: \emph{What is the probability threshold $p=p(n)$ for the Erd\H{o}s--R\'enyi random graph $G_{n,p}$ to contain asymptotically almost surely (a.a.s.) each tree/all trees from a given class $\mathcal F_n$ of trees?} Note that there is a difference between containing ``each tree'' and ``all trees'' (i.e., all trees simultaneously; this is often referred to as \emph{universality}) as the error probabilities for missing individual trees might sum up.

Most research focused on containment of spanning trees, or almost spanning trees. The only well-understood case is when $\mathcal F_n=\{P_{k_n}\}$ is a path. The threshold $p=\frac{(1+o(1))\ln n}{n}$ for appearance of a spanning path (i.e., $k_n=n$) was determined by Koml\'os and Szemer\'edi~\cite{KomSze:HamiltonRandom}, and independently by Bollob\'as~\cite{Boll:HamiltonRandom}. Note that this threshold is the same as the threshold for connectedness. We should also mention a previous result of P\'osa~\cite{Posa:HamiltonRandom} which determined the order of magnitude of the threshold, $p=\Theta(\frac{\ln n}n)$. The heart of P\'osa's proof, the celebrated rotation-extension technique, is an argument about expanding graphs, and indeed many other results about trees in random graphs exploit the expansion properties of $G_{n,p}$ in the first place.

The threshold for the appearance of almost spanning paths in $G_{n,p}$ was determined by Fernandez de la Vega~\cite{Fern:LongRandom} and independently by Ajtai, Koml\'os, and Szemer\'edi~\cite{AKS:LongRandom}. Their results say that a path of length $(1-\epsilon)n$ appears a.a.s.\ in $G_{n,\frac{C}{n}}$ for $C=C(\epsilon)$ sufficiently large. This behavior extends to bounded degree trees. Indeed, Alon, Krivelevich, and Sudakov~\cite{AlKrSu:NearlySpanningTrees} proved that $G_{n,\frac{C}{n}}$ (for a suitable $C=C(\epsilon,\Delta)$) a.a.s.\ contains all trees of order $(1-\epsilon)n$ with maximum degree at most $\Delta$ (the constant $C$ was later improved in~\cite{BCPS:AlmostSpanningRandom}).

Let us now turn to spanning trees in random graphs. It is known~\cite{AlKrSu:NearlySpanningTrees}  that a.a.s.\  $G_{n,\frac{C\ln n}{n}}$ contains a single spanning tree $T$ with bounded maximum degree and linearly many leaves. This result can be reduced to the main result of~\cite{AlKrSu:NearlySpanningTrees} regarding almost spanning trees quite easily. The constant $C$ can be taken $C=1+o(1)$, as was shown recently by Hefetz, Krivelevich, and Szab\'o~\cite{HeKrSz:SpanningRandom}; obviously this is best possible. The same result also applies to trees that contain a path of linear length whose vertices all have degree two. A breakthrough in the area was achieved by Krivelevich~\cite{Kri:SpanningRandom} who gave an upper bound on the threshold $p=p(n,\Delta)$ for embedding a single spanning tree of a given maximum degree $\Delta$. This bound is essentially tight for $\Delta=n^c$, $c\in(0,1)$. Even though the argument in~\cite{Kri:SpanningRandom} is not difficult, it relies on a deep result of 
Johansson, Kahn and Vu~\cite{JoKaVu:Factors} about factors in random graphs. Montgomery~\cite{MontgomeryBoundedDeg} complemented Krivelevich's result obtaining an almost tight upper bound on $p(n,\Delta)$ in the case when $\Delta$ is small. Further, Montgomery~\cite{MontgomeryCombs} achieved an essentially optimal bound for containment of some comb-like graphs.

Regarding universality of random graphs with respect to spanning trees, most of the research focused on the subclass of bounded-degree trees. Let us mention papers~\cite{JoKriSa:Expanders} and~\cite{FerNenPet} which improve the upper-bounds for the probability of containing all trees of maximum degree $\Delta$ (the results are meaningful for $\Delta<n^c$ for some small value of $c$).

\paragraph{Trees in expanders.}
By an expander graph we mean a graph with a large Cheeger constant, i.e., a graph which satisfies a certain isoperimetric property. As indicated above, random graphs are very good expanders, and this is the main motivation for studying tree containment problems in expanders. Another motivation comes from studying the universality phenomenon. Here the goal is to construct sparse graphs which contain all trees from a given class, and expanders are natural candidates for this. The study of sparse tree-universal graphs is a remarkable area by itself which brings challenges both in probabilistic and explicit constructions. For example, Bhatt, Chung, Leighton, and Rosenberg~\cite{BCLR:Universal} give an explicit construction of a graph with only $O_{\Delta}(n)$ edges which contains all $n$-vertex trees with maximum degree at most $\Delta$. The above mentioned paper by
Johannsen, Krivelevich, and 
Samotij~\cite{JoKriSa:Expanders} shows a number of universality results for expanders, too.
 For example, they show universality for the class of graphs with a large Cheeger constant that satisfy a certain connectivity condition.

Friedman and Pippenger~\cite{FP87} extended P\'osa's rotation-extension technique  from paths to trees  and found many applications (e.g.~\cite{HaKo:SizeRamsey,Haxell:TreeEmbeddings,BCPS:AlmostSpanningRandom}). Sudakov and Vondr\'ak~\cite{SudVoTree} use tree-indexed random walks to embed trees in $K_{s,t}$-free graphs (this property implies expansion); a similar approach is employed  by  Benjamini and Schramm~\cite{BeSch:TreeCheeger} in the setting of infinite graphs. Tree-indexed random walks are also used (in conjunction with the Regularity Lemma) in the work of K\"uhn, Mycroft, and Osthus on Sumner's universal tournament conjecture, \cite{KMO:ApproximateSumner,KMO:ProofSumner}.

In our proof of Theorem~\ref{thm:main}, embedding trees in expanders play a crucial role, too. However, our notion of expansion is very different from those studied previously. (Actually, we introduce two, very different, notions  in Definitions~\ref{def:densespot} and~\ref{def:avoiding}.)

\subsection{Overview of the proof of our main result}\label{ssec:OverviewOfProof}
This is a very brief overview of the proof. A more thorough overview is given in~\cite{LKS:overview}.

\medskip

The structure of the proof of our main result (Theorem~\ref{thm:main}) resembles the proof of the 
dense case, Theorem~\ref{thm:PiguetStein}.
We obtain an approximate representation --- called the \emph{sparse decomposition} --- of the host
graph~$G$ from Theorem~\ref{thm:main}.
Then we find a suitable combinatorial structure inside the sparse decomposition.
Finally, we embed a given tree~$T$ into~$G$ using this structure. 

\medskip

Here we explain the key ingredients of the proof in more detail. The input graph $G$ has $\Theta(kn)$ edges. Indeed, an easy counting argument gives that $e(G)\ge kn/4$. On the other hand, we can assume that $e(G)< kn$, as otherwise $G$ contains a subgraph of minimum degree at least $k$, and the assertion of Theorem~\ref{thm:main} follows. Recall that the Szemer\'edi regularity lemma gives an approximation of dense graphs in which $o(n^2)$ edges are neglected. The sparse decomposition introduced here captures all but at most $o(kn)$ edges. The vertex set of $G$ is partitioned into a set of vertices of degree much larger than~$k$ and a set of vertices of degree $O(k)$. 
Further, the induced graph on the second set is split into regular pairs (in the sense of the Szemer\'edi regularity lemma) with clusters of sizes $\Theta(k)$ leading to a cluster graph $\BGblack$, and into two additional 
parts which each have certain (different) expansion properties. The first of these two expanding parts --- called $\Gexp$ --- is a subgraph of $G$ that contains no bipartite subgraphs of a density above a certain threshold density (we call such bipartite subgraphs {\it dense spots}). The second expanding part --- called the \emph{avoiding set $\smallatoms$}--- consists of vertices that lie in many of these dense spots. The vertices of huge degrees, the regular pairs, and the two expanding parts form the sparse decomposition of $G$. 
The key ideas behind obtaining this sparse decomposition are given in~\cite[Section~3]{LKS:overview}, and full details can be found in Section~\ref{sec:class}.
 It is well-known that regular pairs are suitable for embedding small trees. 
In~\cite{cite:LKS-cut3} we work out techniques for embedding small trees in each of the three remaining parts of the sparse decomposition. A nontechnical description of these techniques is given in Section~\ref{sssec:whyavoiding} (for $\smallatoms$) and Section~\ref{sssec:whyGexp} (for $\Gexp$). It is a bit difficult to describe precisely the way the huge degree vertices are utilized. At this moment it suffices to say that it is easy to extend a partial embedding of a $k$-vertex tree from a vertex $u$ mapped to a huge-degree vertex $x$ to the children of $u$. Of course, for such an extension alone, $\deg(x)\ge k-1$ would have been sufficient. So, the fact that the degree of~$x$ is much larger than~$k$ is used (together with other properties) to accommodate these children so that it is possible to continue even with subsequent extensions.

 Tree-embedding results in the dense setting (e.g.\
 Theorem~\ref{thm:PiguetStein}) rely on finding a (connected) matching structure in the
 cluster graph. Indeed, this allows for distributing different parts of the
 tree in the matching edges. In analogy, in the second paper of this series~\cite{cite:LKS-cut1} we find a
 structure based on the sparse decomposition.
 This \emph{rough structure} utilizes all the concepts suitable for embedding trees described above: huge degree vertices, the avoiding set $\smallatoms$, the graph $\Gexp$, and dense regular pairs. Somewhat surprisingly, the dense regular pairs do not come only from $\BGblack$. Let us make this more precise. An initial matching structure is found in~$\BGblack$ and this structure is enhanced using other parts of $G$ to yield further regular pairs, referred to in this context as the \emph{\semiregular matching}. One may ask what the role of $\BGblack$ is. The answer is that either we can take directly a sufficiently large matching in $\BGblack$, or the lack of any such matching in $\BGblack$ gives us information about a compensating enhancement in a form of a \semiregular matching based on other parts of the decomposition. A simplified version of this rough structure is given as Lemma~7 in~\cite{LKS:overview}.

 However, the rough structure is not immediately suitable for embedding $T$, and we shall further refine it in the third paper of this series~\cite{cite:LKS-cut2}. We will show that in the setting of Theorem~\ref{thm:main}, we can always find one of ten \emph{configurations},
 denoted by $\mathbf{(\diamond1)}$--$\mathbf{(\diamond10)}$, in the host graph $G$.
Obtaining these configurations from the rough structure is based on pigeonhole-type arguments such as: if
there are many edges between two sets, and few ``kinds'' of edges, then many of
the edges are of the same kind. The different kinds of edges come from the sparse decomposition (and allow for different kinds of embedding techniques). Just ``homogenizing'' the situation by
restricting to one particular kind is not enough, we also need to
employ certain ``cleaning lemmas''.
A simplest such lemma would be that a graph with many edges contains a subgraph
with a large minimum degree; the latter property evidently being more directly
applicable for a sequential embedding of a tree. The actual cleaning lemmas we
use are complex extensions of this simple idea.

Finally, in~\cite{cite:LKS-cut3}, we show how to embed the tree~$T$. This is done by first establishing some
elementary embedding lemmas for small subtrees, and then combine these
for each of the
cases $\mathbf{(\diamond1)}$--$\mathbf{(\diamond10)}$ to yield an embedding of
the entire tree~$T$.

\begin{figure}[ht!]
\centering 
\includegraphics[scale=0.83]{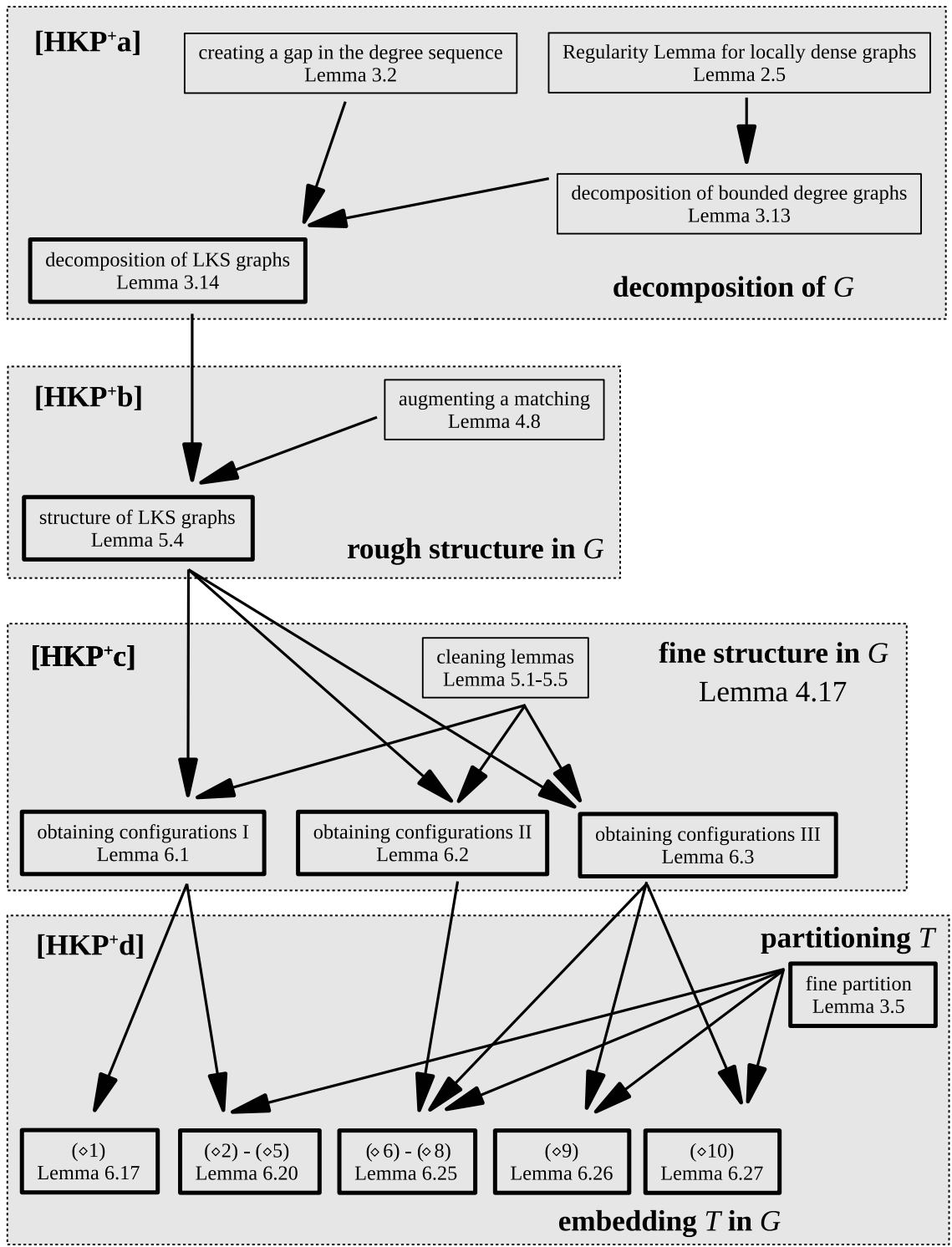}
\caption[Structure of proof of Theorem~\ref{thm:main}]{Structure of the proof of Theorem~\ref{thm:main}, including parts from~\cite{cite:LKS-cut1,cite:LKS-cut2,cite:LKS-cut3}.}
\label{fig:proofstructure}
\end{figure}

A scheme of the proof of Theorem~\ref{thm:main} is given in Figure~\ref{fig:proofstructure}. 

%% file: prelim0.tex
\subsection{General notation}
All graphs considered in this paper are finite, undirected, without multiple edges, and without self-loops. We write \index{mathsymbols}{*VG@$V(G)$}$V(G)$ and \index{mathsymbols}{*EG@$E(G)$}$E(G)$ for the vertex set and edge set of a graph $G$, respectively. Further, \index{mathsymbols}{*VG@$v(G)$}$v(G)=|V(G)|$ is the order of $G$, and \index{mathsymbols}{*EG@$e(G)$}$e(G)=|E(G)|$ is its number of edges. If $X,Y\subset V(G)$ are two, not necessarily disjoint, sets of vertices we write \index{mathsymbols}{*EX@$e(X)$}$e(X)$ for the number of edges induced by $X$, and \index{mathsymbols}{*EXY@$e(X,Y)$}$e(X,Y)$ for the number of ordered pairs $(x,y)\in X\times Y$ such that $xy\in E(G)$. In particular, note that $2e(X)=e(X,X)$.

\index{mathsymbols}{*DEG@$\deg$}\index{mathsymbols}{*MINDEG@$\mindeg$}\index{mathsymbols}{*MAXDEG@$\maxdeg$}
For a graph $G$, a vertex $v\in V(G)$ and a set $U\subset V(G)$, we write
$\deg(v)$ and $\deg(v,U)$ for the degree of $v$, and for the number of
neighbours of $v$ in $U$, respectively. We write $\mindeg(G)$ for the minimum
degree of $G$, $\mindeg(U):=\min\{\deg(u)\::\: u\in U\}$, and
$\mindeg(V_1,V_2)=\min\{\deg(u,V_2)\::\:u\in V_1\}$ for two sets $V_1,V_2\subset
V(G)$. 
Note that for us, the minimum degree of a graph on zero vertices is $\infty$.
Similar notation is used for the maximum degree, denoted by $\maxdeg(G)$.
The neighbourhood of a vertex $v$ is denoted by
\index{mathsymbols}{*N@$\neighbour(v)$}$\neighbour(v)$. We set $\neighbour(U):=\bigcup_{u\in
U}\neighbour(u)$. The symbol $-$ is used
for two graph operations: if $U\subset V(G)$ is a vertex
set then $G-U$ is the subgraph of $G$ induced by the set
$V(G)\setminus U$. If $H\subset G$ is a subgraph of $G$ then the graph
$G-H$ is defined on the vertex set $V(G)$ and corresponds
to deletion of edges of $H$ from $G$.
 Any graph with zero edges is called \index{general}{empty graph}\emph{empty}.

A family $\mathcal A$ of pairwise disjoint subsets of $V(G)$ is an \index{general}{ensemble}\index{mathsymbols}{*ENSEMBLE@ensemble}\emph{$\ell$-ensemble in $G$} if  $|A|\ge \ell$ for each $A\in\mathcal A$. 

The set $\{1,2,\ldots, n\}$ of the first $n$ positive integers is
denoted by \index{mathsymbols}{*@$[n]$}$[n]$. 

Suppose that we have a nonempty set $A$, and $\mathcal X$ and $\mathcal Y$ each partition $A$. Then~\index{mathsymbols}{*@$\boxplus$}$\mathcal X\boxplus\mathcal Y$ denotes the coarsest common refinement of $\mathcal X$ and $\mathcal Y$, i.e., $$\mathcal X\boxplus\mathcal Y:=\{X\cap Y\::\: X\in \mathcal X, Y\in\mathcal Y\}\setminus \{\emptyset\}\;.$$

We frequently employ indexing by many indices. We write
superscript indices in parentheses (such as $a^{(3)}$), as
opposed to notation of powers (such as $a^3$).
We use sometimes subscript to refer to
parameters appearing in a fact/lemma/theorem. For example,
$\alpha_\PARAMETERPASSING{T}{thm:main}$ refers to the parameter $\alpha$ from Theorem~\ref{thm:main}.
We omit rounding symbols when this does not affect the
correctness of the arguments. In overviews we use the symbol~$\ll$ equivalently to the $o(\cdot)$ symbol. 

In Table~\ref{tab:notation} we indicate our notation system (with an outlook to~\cite{cite:LKS-cut1}--\cite{cite:LKS-cut3}).
\begin{table}[h]
\centering
\caption{Specific notation used in the series.}
\label{tab:notation}
\begin{tabular}{r|l}
\hline
lower case Greek letters &  small positive constants ($\ll 1$)\\
                         & $\phi$ reserved for embedding; $\phi:V(T)\rightarrow V(G)$\\
\hline
upper case Greek letters & large positive constants ($\gg 1$)\\
\hline
one-letter bold& sets of clusters \\
\hline
bold (e.g., $\treeclass{k},\LKSgraphs{n}{k}{\eta}$)& classes of graphs\\
\hline
blackboard bold (e.g., $\HugeVertices,\smallatoms,\smallvertices{\eta}{k}{G},\XA$)& distinguished vertex sets, except for\\
& $\NN$ which denotes the set $\{1,2,\ldots\}$\\
\hline 
script  (e.g., $\mathcal A,\mathcal D,\mathcal N$)& families (of vertex sets, ``dense spots'', and regular pairs)\\
\hline
$\class$ (``nabla'')& reserved for ``sparse decomposition''\\
\hline
\end{tabular}
\end{table}

\begin{lemma}\label{lem:subgraphswithlargeminimumdegree}
For all $\ell, n\in\NN$, every $n$-vertex graph $G$ contains a (possibly empty) subgraph $G'$ such that
$\mindeg(G')\ge \ell$ and $e(G')\ge e(G)-(\ell-1) n$.
\end{lemma}
\begin{proof}
We construct the graph $G'$ by sequentially removing vertices of degree less
than $\ell$ from the graph $G$. In each step we remove at most $\ell-1$
edges. Thus the statement follows.
\end{proof}

\input{regul0.tex}

\subsection{LKS graphs}\index{mathsymbols}{*LKSgraphs@$\LKSgraphs{n}{k}{\eta}$}
\label{ssec:LKSgraphs}
Write \index{mathsymbols}{*LKSgraphs@$\LKSgraphs{n}{k}{\eta}$}$\LKSgraphs{n}{k}{\alpha}$
for the class of all $n$-vertex graphs with at least
$(\frac12+\alpha)n$ vertices of degrees at least
$(1+\alpha)k$. With this notation, Conjecture~\ref{conj:LKS} states that every graph in $\LKSgraphs{n}{k}{0}$ contains every tree from $\treeclass{k+1}$.

Given a graph $G$, denote by
\index{mathsymbols}{*S@$\smallvertices{\eta}{k}{G}$}$\smallvertices{\eta}{k}{G}$ the set of those
vertices of $G$ that have degree less than $(1+\eta)k$ and by
\index{mathsymbols}{*L@$\largevertices{\eta}{k}{G}$}$\largevertices{\eta}{k}{G}$ the set of those
vertices of $G$ that have degree at least $(1+\eta)k$.\footnote{``$\mathbb S$'' stands for ``small'', and ``$\mathbb L$'' for ``large''.} Thus the sizes of the sets $\smallvertices{\eta}{k}{G}$ and $\largevertices{\eta}{k}{G}$ are what specifies the membership to $\LKSgraphs{n}{k}{\eta}$.

Define \index{mathsymbols}{*LKSmingraphs@$\LKSmingraphs{n}{k}{\eta}$}$\LKSmingraphs{n}{k}{\eta}$ as the set 
of all graphs $G\in \LKSgraphs{n}{k}{\eta}$ that are  edge-minimal with respect to the membership in $\LKSgraphs{n}{k}{\eta}$.
In order to prove Theorem~\ref{thm:main} it suffices to restrict our attention to graphs from $\LKSmingraphs{n}{k}{\eta}$, and this is why we introduce the class. 
Let us collect some properties of graphs in
$\LKSmingraphs{n}{k}{\eta}$.
\begin{fact}\label{fact:propertiesOfLKSminimalGraphs}
For any graph $G\in\LKSmingraphs{n}{k}{\eta}$ the following is true.
\begin{enumerate}
  \item $\smallvertices{\eta}{k}{G}$ is an independent set.\label{Sisindep}
  \item All the neighbours of every vertex $v\in V(G)$ with $\deg(v)>\lceil(1+\eta)k\rceil$ have degree exactly $\lceil(1+\eta)k\rceil$. \label{neighbinL}
  \item $|\largevertices{\eta}{k}{G}|\le\lceil
  (1/2+\eta)n\rceil+1$.\label{fewlargevs}
\end{enumerate}
\end{fact}
Observe that every edge in a graph $G\in\LKSmingraphs{n}{k}{\eta}$ is incident
to at least one vertex of degree exactly $\lceil(1+\eta)k\rceil$. This gives
the following inequality.
\begin{equation}\label{eq:LKSminimalNotManyEdges}
e(G)\le\lceil(1+\eta)k\rceil
\left|\largevertices{\eta}{k}{G}\right|\leBy{F\ref{fact:propertiesOfLKSminimalGraphs}(\ref{fewlargevs}.)}
\lceil(1+\eta)k\rceil\left(\left\lceil
\left(\frac12+\eta\right)n\right\rceil+1\right)<kn\;.
\end{equation}
(The last inequality is valid under the additional mild
assumption that, say, $\eta<\frac1{20}$ and $n>k>20$. This can be assumed throughout the paper.)

\begin{definition}\label{def:LKSsmall}
Let \index{mathsymbols}{*LKSsmallgraphs@$\LKSsmallgraphs{n}{k}{\eta}$}$\LKSsmallgraphs{n}{k}{\eta}$ be the class of those graphs $G\in\LKSgraphs{n}{k}{\eta}$ for which we have the following three properties:
\begin{enumerate}
   \item All the neighbours of every vertex $v\in V(G)$ with $\deg(v)>\lceil(1+2\eta)k\rceil$ have degrees at most $\lceil(1+2\eta)k\rceil$.\label{def:LKSsmallA}
   \item All the neighbours of every vertex of $\smallvertices{\eta}{k}{G}$
    have degree exactly $\lceil(1+\eta)k\rceil$. \label{def:LKSsmallB}
   \item We have $e(G)\le kn$.\label{def:LKSsmallC}
\end{enumerate}
\end{definition}
Observe that the graphs from $\LKSsmallgraphs{n}{k}{\eta}$ also satisfy~\ref{Sisindep}., and a quantitatively somewhat weaker version of~\ref{neighbinL}.\ of Fact~\ref{fact:propertiesOfLKSminimalGraphs}. This suggests that in some sense $\LKSsmallgraphs{n}{k}{\eta}$ is a good approximation of $\LKSmingraphs{n}{k}{\eta}$.

As said, we will prove Theorem~\ref{thm:main} only for graphs from $\LKSmingraphs{n}{k}{\eta}$. However, it turns out that the structure of $\LKSmingraphs{n}{k}{\eta}$ is too rigid. In particular, $\LKSmingraphs{n}{k}{\eta}$ is not closed under discarding a small amount of edges during our cleaning procedures. This is why the class $\LKSsmallgraphs{n}{k}{\eta}$ comes into play: starting with a graph in $\LKSmingraphs{n}{k}{\eta}$ we perform some initial cleaning and obtain a graph which lies in $\LKSsmallgraphs{n}{k}{\eta/2}$. We then heavily use its structural properties from Definition~\ref{def:LKSsmall} throughout the proof.

%% file: regul0.tex
\subsection{Regular pairs}

In this section we introduce the notion of regular pairs which is central for
Szemer\'edi's regularity lemma and its extension, discussed in Section~\ref{sec:RegL}. We also list some simple properties of regular pairs.

Given a graph $H$ and a pair $(U,W)$ of disjoint
sets $U,W\subset V(H)$ the
\index{general}{density}\index{mathsymbols}{*D@$\density(U,W)$}\emph{density of the pair $(U,W)$} is defined as
$$\density(U,W):=\frac{e(U,W)}{|U||W|}\;.$$
For a given $\varepsilon>0$, a pair $(U,W)$ of disjoint
sets $U,W\subset V(H)$ 
is called an \index{general}{regular pair}\emph{$\epsilon$-regular
pair} if $|\density(U,W)-\density(U',W')|<\epsilon$ for every
$U'\subset U$, $W'\subset W$ with $|U'|\ge \epsilon |U|$, $|W'|\ge
\epsilon |W|$. If the pair $(U,W)$ is not $\epsilon$-regular,
then we call it \index{general}{irregular}\emph{$\epsilon$-irregular}.

We give a useful and well-known property of
regular pairs.
\begin{fact}\label{fact:BigSubpairsInRegularPairs}
Suppose that $(U,W)$ is an $\varepsilon$-regular pair of density
$d$. Let $U'\subset W, W'\subset W$ be sets of vertices with $|U'|\ge
\alpha|U|$, $|W'|\ge \alpha|W|$, where $\alpha>\epsilon$.
Then the pair $(U',W')$ is a $2\varepsilon/\alpha$-regular pair of density at least
$d-\varepsilon$.
\end{fact}

The following fact states a simple relation between the
density of a (not necessarily regular) pair and the densities of its
subpairs.

\begin{fact}\label{fact:CanADensePairConsistOnlyOfSparseSubpairs?}
Let $H=(U,W;E)$ be a bipartite graph of
$\density(U,W)\ge \alpha$. Suppose that the sets $U$ and~$W$ are partitioned into sets $\{U_i\}_{i\in I}$ and
$\{W_j\}_{j\in J}$, respectively. Then at most
$\beta e(H)/\alpha$ edges of $H$ belong to a
pair $(U_i,W_j)$ with $\density(U_i,W_j)\le\beta$.
\end{fact}
\begin{proof}
Trivially, we have
\begin{equation}\label{eq:volumE}
\sum_{i\in I, j\in J}\frac{|U_i||W_j|}{|U||W|}=1\;. 
\end{equation}

Consider a pair $(U_i,W_j)$ of
$\density(U_i,W_j)\le\beta$. Then $$e(U_i,W_j)\le \beta
|U_i||W_j|=\frac{\beta}{\alpha}\frac{|U_i||W_j|}{|U||W|}\alpha|U||W|\le
\frac{\beta}{\alpha}\frac{|U_i||W_j|}{|U||W|}e(U,W)\;.$$
Summing over all such pairs $(U_i,W_j)$ and using~\eqref{eq:volumE} yields the
statement.
\end{proof}


\subsection{Regularizing locally dense graphs}\label{sec:RegL}

The regularity lemma~\cite{Sze78} has proved to be
a powerful tool for attacking graph embedding
problems; see~\cite{KuhnOsthusSurv} 
for a survey. We first state the
lemma in its original form.

\begin{lemma}[Regularity lemma]\label{lem:RL}
For all $\epsilon>0$ and $\ell\in\NN$ there exist $n_0,M\in\NN$
such that for every $n\ge n_0$ the following holds. Let $G$ be an $n$-vertex graph whose
vertex set is pre-partitioned into sets
$V_1,\ldots,V_{\ell'}$, $\ell'\le \ell$. Then there exists
a partition $\{U_0,U_1,\ldots,U_p\}$ of $V(G)$, $\ell<p<M$, with the following properties.
\begin{enumerate}[(1)]
\item For every $i,j\in [p]$ we have $|U_i|=|U_j|$, and  $|U_0|<\epsilon n$.
\item For every
$i\in [p]$ and every $j\in [\ell']$ either $U_i\cap
V_j=\emptyset$ or $U_i\subset
V_j$. 
\item All but at most
$\epsilon p^2$ pairs $(U_i,U_j)$, $i,j\in [p]$, $i\neq j$, are $\epsilon$-regular.
\end{enumerate}
\end{lemma} 
Property~(3) of Lemma~\ref{lem:RL} is often called \emph{$\epsilon$-regularity of the partition $\{U_0,U_1,\ldots,U_p\}$}. For us, it is more convenient to introduce this notion in the bipartite context (in Definition~\ref{def:index}).

We shall use Lemma~\ref{lem:RL}  for auxiliary
purposes only as it is helpful only in the setting of dense
graphs (i.e., graphs which have $n$ vertices and
$\Omega(n^2)$ edges). This is not necessarily the
case in Theorem~\ref{thm:main}. For this reason, we  give a version of the 
regularity lemma --- Lemma~\ref{lem:sparseRL} below --- which 
allows us to regularize even sparse graphs.

\begin{figure}[t]
\centering 
\includegraphics[scale=0.8]{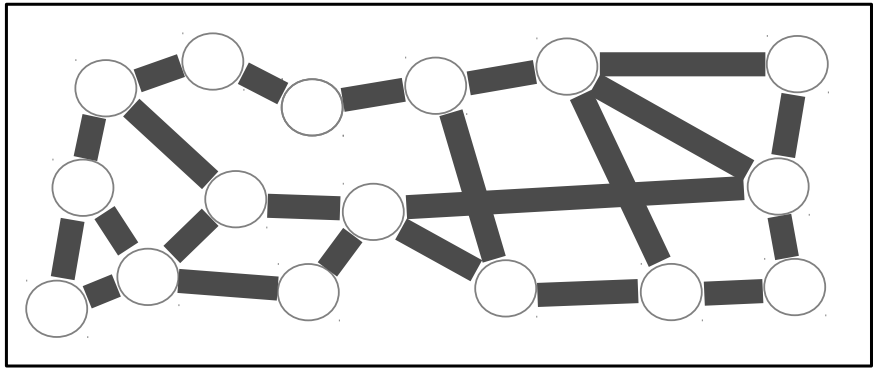}
\caption[Locally dense graph]{A locally dense graph as in Lemma~\ref{lem:sparseRL}. The sets $W_1,\ldots,W_\ell$ are depicted with grey circles. Even though there is a large number of them, each $W_i$ is linked to only boundedly many other $W_j$'s (at most four, in this example). Lemma~\ref{lem:sparseRL} allows us to regularize all the bipartite graphs using the same system of partitions of the sets $W_i$.}\label{fig:locallydensegraph}
\end{figure}

More precisely, suppose that we have an $n$-vertex graph $H$ whose edges lie in bipartite graphs $H[W_i,W_j]$, where $\{W_1,\ldots,W_\ell\}$ is an ensemble of sets of individual sizes $\Theta(k)$. 
Although $\ell$ may be unbounded, for a fixed $i\in[\ell]$ there are only a bounded number (independent of $k$), say $m$, of indices $j\in[\ell]$ such that  $H[W_i,W_j]$ is non-empty. See Figure~\ref{fig:locallydensegraph} for an example.
Lemma~\ref{lem:sparseRL} then allows us to regularize (in the sense of the regularity lemma, Lemma~\ref{lem:RL}) all the bipartite graphs $G[W_i,W_j]$ using the same partition $\{W_i^{(0)}\dcup W_i^{(1)}\dcup\ldots\dcup W_i^{(p_i)}=W_i\}_{i=1}^\ell$. Note that as $|W_i|=\Theta(k)$ for all $i\in[\ell]$ then $H$ has at most
$$\Theta(k^2)\cdot m\cdot \ell\le \Theta(k^2)\cdot m\cdot \frac n{\Theta(k)}=\Theta(kn)$$
edges. Thus, when $k\ll n$, this is a regularization of a sparse graph. This ``sparse regularity lemma'' is very different to that of Kohayakawa and R\"odl (see e.g.~\cite{Kohayakawa97Szemeredi}). Indeed, the Kohayakawa--R\"odl regularity lemma only deals with graphs which have no local condensation of edges, e.g., subgraphs of random graphs.\footnote{There is a recent refinement of the Kohayakawa--R\"odl regularity lemma, due to Scott~\cite{Scott:RL}. Scott's regularity lemma gets around the no-condensation condition, which proves helpful in some situations, e.g.~\cite{AllKeeSudVer:TuranBipartite}; still, the main features remain.} Consequently, the resulting regular pairs are of density $o(1)$. In contrast, Lemma~\ref{lem:sparseRL} provides us with regular pairs of density $\Theta(1)$, but, on the other hand, is useful only for graphs which are locally dense.

\begin{lemma}[Regularity lemma for locally dense graphs]\label{lem:sparseRL}
For all $m,z\in \mathbb{N}$ and $\epsilon>0$  there
exists $q_\mathrm{MAXCL}\in\NN$ such that the following is true. Suppose $H$ and $F$ are two graphs, $V(F)=[\ell]$ for some $\ell\in\mathbb N$, and 
$\maxdeg(F)\le m$. Suppose that $\mathcal Z=\{Z_1,\ldots,Z_z\}$ is a partition of $V(H)$. Let $\{W_1,\ldots, W_\ell\}$ 
be a $q_\mathrm{MAXCL}$-ensemble in $H$, such that for all
$i,j\in[\ell]$ we have
\begin{equation}
\label{lem:sparseRL(item)samesize} 
2|W_i|\ge |W_j|\;.
\end{equation}
Then for each $i\in [\ell]$ there exists a partition
$W_i^{(0)},W_i^{(1)},\ldots,W_i^{(p_i)}$ of the set
$W_i$  such that for all
$i,j\in[\ell]$ we have
\begin{enumerate}[(a)]
\item\label{item:RL:a} $1/\epsilon\le p_i\le q_\mathrm{MAXCL}$,
\item $|W_i^{(i')}|=|W_j^{(j')}|$ for each $i'\in [p_i]$,
$j'\in[p_j]$, 
\item for each $i'\in [p_i]$ there exists $x\in[z]$ such that $W_i^{(i')}\subset Z_x$,
\item\label{item:RL:garbage} $\sum_i|W_i^{(0)}|<\epsilon \sum_i|W_i|$, and\label{RLgarbage}
\item\label{item:RL:d} at most $\epsilon\left|\mathcal{Y}\right|$ pairs
$\left(W_i^{(i')},W_j^{(j')}\right)\in \mathcal{Y}$ form an
$\epsilon$-irregular pair in $H$, where
$$\mathcal{Y}:=\left\{\left(W_i^{(i')},W_j^{(j')}\right)\::\:ij\in
E(F), i'\in [p_i], j'\in [p_j]\right\}\;.$$
\end{enumerate}
\end{lemma}
We use Lemma~\ref{lem:sparseRL} in the proof of Lemma~\ref{lem:decompositionIntoBlackandExpanding}. Lemma~\ref{lem:decompositionIntoBlackandExpanding} is in turn the main tool in the proof of our main structural decomposition of the graph $G_\PARAMETERPASSING{T}{thm:main}$, Lemma~\ref{lem:LKSsparseClass}. In the proof of Lemma~\ref{lem:LKSsparseClass} we decompose the input graph into several parts with very different properties, and one of these parts is a locally dense graph which can be then regularized by Lemma~\ref{lem:decompositionIntoBlackandExpanding}.
A similar regularity lemma is used
in~\cite{AKSS07+}. 
\bigskip 

The proof of Lemma~\ref{lem:sparseRL} is similar to the proof of the standard
regularity lemma (Lemma~\ref{lem:RL}), as given for example in~\cite{Sze78}. The key notion is that of the index (a.k.a.\ the mean square density) which we recall now. For us, it is convenient to work in the category of bipartite graphs.
\begin{definition}\label{def:index}
Suppose that $\mathcal X=\{X_0,X_1,\ldots,X_\ell\}$ and $\mathcal Y=\{Y_0,Y_1,\ldots,Y_p\}$ are partitions of a set $X$, and of $Y$ with distinctive sets $X_0$ and $Y_0$ which we call \index{general}{garbage cluster}\emph{garbage clusters}. We use the symbol $\circ$ to indicate a new partition in which the garbage cluster is broken into singletons, e.g., $\mathcal X^\circ=\{X_1,\ldots,X_\ell\}\cup\{ \{x\}:x\in X_0\}$. We say that \index{general}{refine up to garbage cluster}\emph{$\mathcal X$ refines $\mathcal Y$ up to the garbage cluster} if $\mathcal X^\circ$ refines $\mathcal Y^\circ$.

Suppose that $G=(A,B;E)$ is a bipartite graph. Let $\mathcal A=\{A_0,A_1,\ldots,A_s\}$ and $\mathcal B=\{B_0,B_1,\ldots,B_t\}$ be partitions of $A$ and $B$, with garbage clusters $A_0$ and $B_0$. We say that the pair $(\mathcal A,\mathcal B)$ is an \index{general}{regular partition}\emph{$\epsilon$-regular partition} of $G$ if at most $\epsilon st$ pairs $(A_i,B_j)$, $i\in[s],j\in[t]$ are irregular. Otherwise, $(\mathcal A,\mathcal B)$ is \index{general}{irregular partition}\emph{$\epsilon$-irregular}.

We then define the \index{general}{index}\emph{index} of $(\mathcal A,\mathcal B)$ by
\index{mathsymbols}{*IND@$\mathrm{ind}(\mathcal A,\mathcal B)$}
$$\mathrm{ind}(\mathcal A,\mathcal B)=\frac1{(|A|+|B|)^2}\cdot \sum_{X\in \mathcal A^\circ, Y\in \mathcal B^\circ}\frac{e(X,Y)^2}{|X||Y|}\;.$$
\end{definition}
Clearly, $\mathrm{ind}(\mathcal A,\mathcal B)\in[0,1]$. Here is another fundamental property of the index.
\begin{fact}[Bipartite version of Lemma~7.4.2 in~\cite{Die05}]\label{fact:index}
Suppose that $G=(A,B;E)$ is a bipartite graph. Let $\mathcal A$ and $\mathcal A'$ be partitions of $A$ with given garbage clusters. Let $\mathcal B$ and $\mathcal B'$ be partitions of $B$ with given garbage clusters. Suppose that $\mathcal A'$ refines $\mathcal A$ and $\mathcal B'$ refines $\mathcal B$ up to garbage clusters. Then $\mathrm{ind}(\mathcal A',\mathcal B')\ge \mathrm{ind}(\mathcal A,\mathcal B)$.
\end{fact}


\begin{lemma}[Index Pumping Lemma; bipartite version of Lemma~7.4.4.\ in~\cite{Die05}]\label{lem:IndexPump}
Let $\epsilon\in (0,\frac14)$ and $p,q\in\NN$. Let $G$ be a bipartite graph $G=(A,B;E)$, with $\frac{|A|}2\le |B|\le 2|A|$. Suppose that $\mathcal A$ and $\mathcal B$ are partitions of vertex sets $A$ and $B$ with distinctive garbage clusters $A_0$ and $B_0$. Suppose further that
\begin{enumerate}[(a)]
	\item $p \le|\mathcal A|,|\mathcal B|\le q$, \item $|A_0|<\epsilon |A|$, $|B_0|<\epsilon |B|$, and
	\item all the sets in $\mathcal A\cup \mathcal B\setminus\{A_0,B_0\}$ have the same size.
\end{enumerate}
If $(\mathcal A,\mathcal B)$ is not an $\epsilon$-regular partition of $G$ then there exist partitions $\mathcal A'$ and $\mathcal B'$ of $A$ and $B$ with garbage clusters $A_0'$ and $B_0'$ such that 
\begin{enumerate}[(i)]
	\item $p+1\le |\mathcal A'|,|\mathcal B'|\le 2q16^q$, and
	\item $|A_0'|\le |A_0|+\frac{|A|}{2^p}$, $|B_0'|\le |B_0|+\frac{|B|}{2^p}$,
	\item
	all the sets in $\mathcal A'\cup \mathcal B'\setminus\{A_0',B_0'\}$ have the same size, 
	\item the partitions $\mathcal A'$ and $\mathcal B'$ refine $\mathcal A$ and $\mathcal B$ up to garbage clusters, and
	\item $\mathrm{ind}(\mathcal A',\mathcal B')\ge \mathrm{ind}(\mathcal A,\mathcal B)+\tfrac{\epsilon^5}{3691}$.
\end{enumerate}   
\end{lemma}
We note that by stating a version for bipartite graphs we had to adjust numerical values compared to~\cite{Die05}. Recall that the proof of Lemma~\ref{lem:IndexPump} has two independent steps: first the partitions $\mathcal A$ and $\mathcal B$ are suitably refined (so that the index increases) and then these new partitions are further refined (up to garbage clusters) so that the non-garbage sets have the same size. The latter step does not decrease the index by Fact~\ref{fact:index} but may possibly increase the sizes of the garbage clusters. Thus, we can state a version of Lemma~\ref{lem:IndexPump} in which refinements are performed simultaneously on a number of bipartite graphs (referred to as $(G_i)_i$ in the corollary below), and in addition refines further partitions (referred to as $(C_j)_j$ below) on which no regularization is imposed.
\begin{cor}\label{cor:IndexPumpSimult}
Let $\epsilon\in (0,\frac14)$ and $p,q\in\NN$. Let $G_i$ $i\in I$ be bipartite graphs $G_i=(A_i,B_i;E_i)$. Let $C_j$, $j\in J$ be sets of vertices. Suppose that all the sets $A_i$, $B_i$, and $C_j$ are mutually disjoint. Suppose further that for each $i\in I$ and $j\in J$, $\max\{\frac{|A_i|}2,\frac{|B_i|}2\}\le |C_j|\le \min\{2|A_i|,2|B_i|\}$.  Suppose that $\mathcal A_i$ and $\mathcal B_i$ are partitions of $A_i$ and $B_i$ with garbage clusters $A_{0i}$ and $B_{0i}$, and that $\mathcal C_j$ are partitions of $C_j$ with garbage clusters $C_{0j}$. Suppose further that for each $i\in I$, and $j\in J$,
\begin{enumerate}[(a)]
	\item\label{cor:boundABass} $p \le|\mathcal A_i|,|\mathcal B_i|,|\mathcal C_j|\le q$, \item $|A_{0i}|<\epsilon |A_i|$, $|B_{0i}|<\epsilon |B_i|$, $|C_{0j}|<\epsilon |C_j|$, and
	\item all the sets in $\bigcup_{m\in I}(\mathcal A_m\cup \mathcal B_m\setminus\{A_{0m},B_{0m}\})\cup \bigcup_{n\in J}(\mathcal C_n\setminus\{C_{0n}\})$ have the same size.
\end{enumerate}
If all partitions $(\mathcal A_i,\mathcal B_i)$, $i\in I$ are $\epsilon$-irregular then there exist partitions $\mathcal A_i'$ and $\mathcal B_i'$ of $A$ and $B$ with garbage clusters $A_{0i}'$ and $B_{0i}'$, and partitions $\mathcal C'_j$ of $C_j$ with garbage clusters $C'_{0j}$ such that for each $i\in I$ and $j \in J$,
\begin{enumerate}[(i)]
	\item\label{cor:boundABconc} $p+1\le |\mathcal A'_i|,|\mathcal B'_i|,|\mathcal C'_j|\le 2q\cdot 16^q$, and
	\item $|A_{0i}'|\le |A_{0i}|+\frac{|A_i|}{2^p}$, $|B_{0i}'|\le |B_{0i}|+\frac{|B_i|}{2^p}$, and $|C_{0j}'|\le |C_{0j}|+\frac{|C_j|}{2^p}$, 
	\item
	all the sets in $\bigcup_{m\in I}\mathcal A_m'\cup \mathcal B_m'\setminus\{A_{0m}',B_{0m}'\}\cup \bigcup_{n\in J}\mathcal C_n'\setminus\{C_{0n}'\}$ have the same size, 
	\item the partitions $\mathcal A_i'$, $\mathcal B_i'$ and $\mathcal C_j'$ refine $\mathcal A_i$, $\mathcal B_i$ and $\mathcal C_j$ up to garbage clusters, and
	\item $\mathrm{ind}(\mathcal A_i',\mathcal B_i')\ge \mathrm{ind}(\mathcal A_i,\mathcal B_i)+\tfrac{\epsilon^5}{3691}$.
\end{enumerate}   
\end{cor}

We are now in a position when we can prove Lemma~\ref{lem:sparseRL}. But before, let us describe how a more naive
approach fails. For each edge $ij\in E(F)$ consider a regularization of the
bipartite graph $H[W_i,W_j]$,  let
$\{U^{(i')}_{i,j}\}_{i'\in[q_{i,j}]}$ be the partition of $W_i$ into clusters, and let
$\{U^{(j')}_{j,i}\}_{j'\in[q_{j,i}]}$ be the partition of $W_j$ into
clusters such that almost all pairs $(U^{(i')}_{i,j},U^{(j')}_{j,i})\subseteq (W_i,W_j)$ form an
$\epsilon'$-regular pair (for some $\epsilon'$ of our taste). We would now be
done if the partition $\{U^{(i')}_{i,j}\}_{i'\in[q_{i,j}]}$ of $W_i$ was
independent of the choice of the edge $ij$. This however need not be the case.
The natural next step would therefore be to consider the common refinement
$$\underset{j:ij\in E(F)}{\bigboxplus} \big\{U^{(i')_{i,j}}\big\}_{i'\in
[q_{ij}]}$$
of all the obtained partitions of $W_i$. The pairs obtained in this way lack
however any regularity properties as they are too small. Indeed, it is a notorious drawback of the
regularity lemma that the number of clusters in the partition is enormous as a
function of the regularity parameter. In our setting, this means that
$q_{i,j}\gg\frac1{\epsilon'}$. Thus a typical cluster $U^{(i'_1)}_{i,j_1}$
occupies on average only a $\frac1{q_{i,j_1}}$-fraction of the cluster
$U^{(i'_2)}_{i,j_2}$, and thus already the set
$U^{(i'_1)}_{i,j_1}\cap U^{(i'_2)}_{i,j_2}\subset U^{(i'_2)}_{i,j_2}$ is not
substantial (in the sense of the regularity). The same issue arises when regularizing multicoloured graphs
(cf.~\cite[Theorem~1.18]{KS96}). The solution is to impel
the regularizations to happen in a synchronized way.

\begin{proof}[Proof of Lemma~\ref{lem:sparseRL}] 
%

Without loss of generality, assume that $\epsilon<1$. Set $\tilde\epsilon=\epsilon/8$. The number $q_\mathrm{MAXCL}$ can be taken by considering the function $q\mapsto 2q\cdot 16^q$ with initial value $\lceil\frac{4z}{\tilde\epsilon}\rceil$ and iterating it  $\lceil\frac{3691(m+1)}{\tilde\epsilon^6}\rceil$-many times.

For each $i\in [\ell]$ consider an arbitrary initial partition $\mathcal W_i=\{W_i^{(0)},W_i^{(1)},\ldots,W_i^{(p_i)}\}$ of $W_i$ such that for the garbage cluster we have $|W_i^{(0)}|\le \tilde\epsilon|W_i|$, all the non-garbage clusters $W_i^{(i')}$ are disjoint subsets of some set $Z_r$, $r\in [z]$. Further, we make all the non-garbage clusters (coming from all the sets $W_i$) have the same size. It is clear that this can be achieved, and that we can further impose that 
\begin{equation}\label{eq:pi}
1+\frac1\epsilon\le p_i\le \frac{4z}{\epsilon}\;.
\end{equation}

By Vizing's Theorem we can cover the
edges of $F$ by non-empty disjoint matchings
$M_1,\ldots,M_{Q}$, $Q\le m+1$.
For each
$q\in[Q]$ we shall introduce a variable $\mathrm{ind}_q$,
\begin{align*}
\mathrm{ind}_q&=\frac1{|M_q|}\sum_{xy\in M_q}
\mathrm{ind}(\mathcal{W}_x,\mathcal{W}_y)\;.
\end{align*}

We shall now keep refining the partitions $(\mathcal W_z)_{z\in V(F)}$ in steps $\ell=1,2,\ldots$, as follows. Suppose that for some $q\in [Q]$, for the matching $$M_q'=\left\{xy\in M_q\: :\: \mbox{the partition $(\mathcal W_x,\mathcal W_y)$ of the bipartite graph $H[W_x,W_y]$ is $\tilde\epsilon$-irregular}\right\}$$ we have $|M_q'|\ge \tilde\epsilon |M_q|$. 
We apply Corollary~\ref{cor:IndexPumpSimult} with the following setting. The partitions $\{(\mathcal W_x,\mathcal W_y)\}_{xy\in M_q'}$ play the role of $(\mathcal A_i,\mathcal B_i)$'s and the partitions $\{\mathcal W_x\}_{x\in V(F)\setminus V(M_q')}$ play the role of $\mathcal C_j$'s. Further, in step $\ell$ we set $p_\PARAMETERPASSING{C}{cor:IndexPumpSimult}=\frac {1}{\varepsilon}+\ell$. Corollary~\ref{cor:IndexPumpSimult} says that for the modified partitions (which we still denote the same), the index along each edge of $M_q'$ increased by at least $\frac{\tilde\epsilon^5}{3691}$. 
Combined with Fact~\ref{fact:index} and with the fact that $|M_q'|\ge \tilde\epsilon |M_q|$, we get that $\mathrm{ind}_q$ increased by at least $\frac{\tilde\epsilon^6}{3691}$, and none other index $\mathrm{ind}_{t}$ decreased. Observe also that the lower-bound in Corollary~\ref{cor:IndexPumpSimult}\eqref{cor:boundABconc} makes it possible to apply Corollary~\ref{cor:IndexPumpSimult} with $p_\PARAMETERPASSING{C}{cor:IndexPumpSimult}$ increased by one in a next step.

Since $\sum_{t=1}^Q \mathrm{ind}_{t}\le |Q|\le m+1$, we conclude that after at most $\frac{3691(m+1)}{\tilde\epsilon^6}$ steps, for each $q\in [Q]$ the number of bipartite graphs $H[W_x,W_y]$, $xy\in M_q$ that are partitioned $\tilde\epsilon$-irregularly is less than $\tilde\epsilon|M_q|$. In particular, among all the bipartite graphs $(H[W_x,W_y])_{xy\in E(F)}$, at most $\tilde{\epsilon}\cdot e(F)$ are partitioned $\tilde\epsilon$-irregularly.  We claim that this system of partitions satisfies Properties~\eqref{item:RL:a}--\eqref{item:RL:d} as in the statement of the lemma.

Each irregular pair counted in \eqref{item:RL:d} is a pair contained either in $\tilde{\epsilon}$-regularly partitioned or  an $\tilde{\epsilon}$-irregularly partitioned bipartite graph $H[W_x,W_y]$.
It follows from above that the number of irregular pairs of each of these two types is upper-bounded by $\frac{\epsilon}{2}|\mathcal Y|$. For the bound~\eqref{item:RL:garbage}, recall that initially we had $|W_i^{(0)}|\le \frac{\epsilon}2|W_i|$. 
During each application of Corollary~\ref{cor:IndexPumpSimult}, the garbage sets $W_i^{(0)}$ could have grown by at most $\frac{|W_i|}{2^p}$, for  $p=\frac1\epsilon+1,\frac1\epsilon+2,\ldots$. Thus, at the end of the process, we have
$|W_i^{(0)}|\le (\frac{\epsilon}2+\sum_{p>\frac1\epsilon}2^{-p})|W_i|\le {\epsilon}|W_i|$, as needed. The other assertions of the lemma are clear.
\end{proof}

Usually after applying the regularity lemma to some graph $G$, one bounds the
number of  edges which correspond to irregular pairs, to regular, but sparse
pairs, or are incident with the exceptional sets $U_0$. We shall do the same
for the setting of Lemma~\ref{lem:sparseRL}.

\begin{lemma}\label{notmuchlost}
In the situation of Lemma~\ref{lem:sparseRL}, suppose that  
$\maxdeg(H)\leq \Omega k$ and $e(H)\le kn$, and that each edge $xy\in E(H)$ is
captured by some edge $ij\in E(F)$, i.e., $x\in W_i$, $y\in W_j$. Moreover
suppose that
\begin{equation}\label{Fdensi}
\text{$\density(W_i,W_j)\geq\gamma$ if $ij\in E(F)$.}
\end{equation}
Then all but at most $(\frac{4\epsilon}{\gamma}+\epsilon\Omega+\gamma)nk$  edges of $H$ belong to regular pairs $(W^{(i)}_{i'},W^{(j)}_{j'})$, $i,j\neq 0$,  of density at least
$\gamma^2$.
\end{lemma}

\begin{proof}Set $w:=\min\{|W_i|:i\in V(F)\}$.
By~\eqref{Fdensi}, each edge of $F$ represents at
least $\gamma w^2$ edges of $H$. Since $e(H)\le kn$ it follows that $e(F)\le
kn/(\gamma w^2)$. Thus,
by the assumption~\eqref{lem:sparseRL(item)samesize}, $\sum_{AB\in E(F)}|A||B|\le e(F)(2w)^2\le \frac{4kn}{\gamma}$.
Using~\eqref{item:RL:d} of Lemma~\ref{lem:sparseRL} we get that the number of
edges of $H$ contained in $\epsilon$-irregular pairs from $\mathcal
Y$ is at most
\begin{equation}\label{eq:boundTotalIrregII}
\frac{4\epsilon nk}{\gamma}\;.
\end{equation}

Write $E_1$ for the set of edges of $H$ which are
incident with a vertex in $\bigcup_{i\in[\ell]} W_i^{(0)}$. Then
by~\eqref{RLgarbage} of Lemma~\ref{lem:sparseRL}, and since
$\maxdeg(H)\leq\Omega k$,
\begin{equation}\label{eq_incidentToGarbageCluster}
|E_1|\le \epsilon\Omega nk\;.
 \end{equation}
 
Let $E_2$ be the set of those edges of $H$ which belong to $\epsilon$-regular pairs $(W^{(i')}_{i},W^{(j')}_{j})$ with
$ij\in E(F), i'\in[p_i],j'\in[p_j]$ of density at most
$\gamma^2$. We claim that
\begin{equation}\label{eq_edgesinsparseregular}
|E_2|\le\gamma kn\;.
 \end{equation}
Indeed, because of~\eqref{Fdensi} and by
Fact~\ref{fact:CanADensePairConsistOnlyOfSparseSubpairs?}
(with
$\alpha_\PARAMETERPASSING{F}{fact:CanADensePairConsistOnlyOfSparseSubpairs?}:=\gamma$
and
$\beta_\PARAMETERPASSING{F}{fact:CanADensePairConsistOnlyOfSparseSubpairs?}:=\gamma^2$),
for each $ij\in E(F)$ there are at most $\gamma e_H(W_i,W_j)$ edges contained in the
bipartite graphs $H[W^{(i')}_{i},W^{(j')}_{j}]$, 
$i'\in[p_i],j'\in[p_j]$, with $\density_H(W^{(i')}_{i},W^{(j')}_{j})\le \gamma^2$. Since $\sum_{ij\in E(F)}e_H(W_i,W_j)\le kn$, the validity of~\eqref{eq_edgesinsparseregular} follows. Combining \eqref{eq:boundTotalIrregII},~\eqref{eq_incidentToGarbageCluster}, and~\eqref{eq_edgesinsparseregular} we finish the proof.
\end{proof}

%% file: Class0.tex
In this section, we work out a structural decomposition of a possibly
sparse graph which is suitable for embedding trees. Our motivation comes from
the success of the regularity method in the  setting of dense graphs
(see~\cite{KuhnOsthusSurv}). The main technical result of this section, the
``decomposition lemma'', Lemma~\ref{lem:decompositionIntoBlackandExpanding},
provides such a decomposition. Roughly speaking, each graph of a moderate
maximum degree can be decomposed into regular pairs, and two different expanding
parts.

We then combine Lemma~\ref{lem:decompositionIntoBlackandExpanding} with a lemma
on creating a gap in the degree sequence (Lemma~\ref{prop:gap}) to get a decomposition lemma
for graphs from~$\LKSgraphs{n}{k}{\eta}$, Lemma~\ref{lem:LKSsparseClass}.
Lemma~\ref{lem:LKSsparseClass} asserts that each graph
from~$\LKSgraphs{n}{k}{\eta}$ can be decomposed into vertices of degree much
larger than $k$, regular pairs, and expanding parts. 
Further we give a non-LKS-specific version of Lemma~\ref{lem:LKSsparseClass} in Lemma~\ref{lem:genericBD},
which asserts that \emph{each} graph with average degree bigger than an absolute constant has a sparse decomposition. Such a decomposition lemma was used by Ajtai, Koml\'os, Simonovits and Szemer\'edi in their work on the Erd\H os--S\'os conjecture and we expect that it will find applications in other tree embedding problems, and possibly elsewhere.

\subsection{Creating a gap in the degree sequence}\label{ssec:class-gap}
The goal of this section is to show that any graph
$G\in\LKSmingraphs{n}{k}{\eta}$ has a subgraph  $G'\in
\LKSsmallgraphs{n}{k}{\eta/2}$ which has a gap in its degree
sequence. Note that $G'$ then contains almost all the edges of $G$. This is formulated in Lemma~\ref{prop:gap}. Before stating and proving it, we illustrate our proof technique on a simpler version of Lemma~\ref{prop:gap} that applies to all graphs. This simpler lemma will not be used except in the proof of Lemma~\ref{lem:genericBD} which serves also for illustration only.
\begin{lemma}\label{prop:gapIllustration}
Let $(\Omega_i)_{i\in\NN}$ be a sequence of positive numbers with
$\frac{\Omega_j}{\Omega_{j+1}}\leq \frac\eta{2}$ for all $j\in\NN$. Let $G$ be a graph of order $n$ with average degree $k$. Then there is an index $i^{*}\leq \frac4\eta$ and a spanning subgraph $G'\subseteq G$ with $e(G')\ge e(G)-\eta kn$ and with the property that $G'$ contains no vertex with degree in the interval $[\Omega_i k,\Omega_{i+1}k)$.
\end{lemma}
\begin{proof}
Set $R:=\lfloor 4\eta^{-1}\rfloor$. 
For $i\in [R]$ and any graph $H\subseteq G$ define the sets $X_i(H):=\{v\in
V(H)\::\: \deg_H(v)\in[\Omega_i k, \Omega_{i+1}k)\}$ and for $i=R+1$ set
$X_i(H):=\{v\in V(H)\::\: \deg_H(v)\in[\Omega_i k, \infty)\}$. As
\begin{equation*}
\sum_{i\in[R]}\quad\sum_{v\in X_{i}(G)\cup X_{i+1}(G)}\deg(v)\le 4e(G)\;,
\end{equation*}
by averaging  we find an index $i^*\in [R]$ such that
\begin{equation*}
\sum_{v\in X_{i^*}(G)\cup X_{i^*+1}(G)}\deg(v)\le \frac{4e(G)}{R}=\frac{2kn}{R}.
\end{equation*}
Let $G_0\subset G$ be obtained from $G$ by deleting all the edges incident with $X_{i^*}(G)\cup X_{i^*+1}(G)$. In particular,
\begin{equation}\label{eq:GG0}
e(G_0)\ge e(G)-\eta kn/2\;.
\end{equation}
We continue successively deleting edges as follows. If in some step $j=1,2,\ldots$ the set $X_{i^*}(G_{j-1})$ is non-empty, we take an arbitrary vertex $v_j\in X_{i^*}(G_{j-1})$ and obtain a new graph $G_j$ from $G_{j-1}$ by deleting all the (at most $\Omega_{i^*+1} k$ many) edges incident with $v_j$. Obviously, this procedure will terminate eventually. Let $G'$ denote the final graph. Clearly, $G'$ has the desired gap in the degree sequence. It therefore suffices to upper bound $e(G)-e(G')$.

Observe that for any vertex $v_j$ above we have $v_j\in \bigcup_{i=i^*+2}^{R+1} X_i(G)$. Thus,
\begin{equation*}
e(G')-e(G_0)\le 
\Omega_{i^*+1} k\left|\bigcup_{i=i^*+2}^{R+1} X_i(G)\right| 
\le \Omega_{i^*+1} k \cdot\frac{2e(G)}{\Omega_{i^*+2}k}\le \frac{\eta kn}{2}\;.
\end{equation*}
Combined with~\eqref{eq:GG0} we get the statement.
\end{proof}

\begin{proposition}\label{prop:gap}
Let $\eta\in(0,1)$, $G\in\LKSmingraphs{n}{k}{\eta}$ and let $(\Omega_i)_{i\in\NN}$
be a sequence of positive numbers with $\Omega_1>2$ and
$\Omega_j/\Omega_{j+1}\leq \eta^2/100$ for all $j\in\NN$. Then there exist an index $i^{*}\leq 100\eta^{-2}$ and  a subgraph $G'\subseteq G$ such that 
\begin{enumerate}[(i)]
\item\label{item:(prop:gap)stillLKS}$G'\in \LKSsmallgraphs{n}{k}{\eta/2}$, and
\item\label{item:(prop:gap)gap} no vertex $v\in V(G')$ has degree
$\deg_{G'}(v)\in[\Omega_{i^*}k,\Omega_{i^{*}+1}k)$.
\end{enumerate}
\end{proposition}

\begin{proof}
Set $R:=\lfloor 100 \eta^{-2}\rfloor$. 
For $i\in [R]$ and any graph $H\subseteq G$ define the sets $X_i(H):=\{v\in
V(H)\::\: \deg_H(v)\in[\Omega_i k, \Omega_{i+1}k)\}$ and for $i=R+1$ set
$X_i(H):=\{v\in V(H)\::\: \deg_H(v)\in[\Omega_i k, \infty)\}$. As
\begin{equation*}
\sum_{i\in[R]}\quad\sum_{v\in X_{i}(G)\cup X_{i+1}(G)}\deg(v)\le 4e(G)\;,
\end{equation*}
by averaging  we find an index $i^*\in [R]$ such that
\begin{equation}\label{eq:sparselyconnectedpair}
\sum_{v\in X_{i^*}(G)\cup X_{i^*+1}(G)}\deg(v)\le \frac{4e(G)}{R}.
\end{equation}

Let $E_0$ be the set of all the edges incident with $X_{i^*}(G)\cup X_{i^*+1}(G)$. Now, starting with $G_0:=G-E_0$, inductively define graphs $G_j\subsetneq G_{j-1}$ for $j\geq 1$ using any of the following two types of edge deletions:
\begin{enumerate}
\item[(T1)] If there is a vertex $v_j\in X_{i^*}(G_{j-1})$ then we choose an edge $e_j$ 
incident with $v_j$, and set $G_j:=G_{j-1}-e_j$.
\item[(T2)] If there is an edge $e_j=u_jv_j$ of $G_{j-1}$ with $u_j\in \smallvertices{\eta/2}{k}{G_{j-1}}$ and $v_j\in \bigcup_{i=i^*+1}^{R+1}X_{i}(G_{j-1})$ then set $G_j:=G_{j-1}-e_j$. 
\end{enumerate}
Since we keep deleting edges, the procedure stops at some point, say at step $j^*$, when neither of (T1), (T2) is applicable. Note that the resulting graph $G_{j^*}$ already has Property~\eqref{item:(prop:gap)gap}.

 Let $E_1 \subset E(G)$ be the set of those edges deleted by applying (T1). We shall estimate the size of $E_1$. First, observe that
\begin{equation*}
\left|\bigcup_{i=i^*+2}^{R+1} X_i(G)\right| \le \frac{2e(G)}{\Omega_{i^*+2}k}\;.
\end{equation*}
Moreover, each vertex of $\bigcup_{i=i^*+2}^{R+1} X_i(G)$ appears at most $ (\Omega_{i^*+1}-\Omega_{i^*})k < \Omega_{i^*+1}k$ times as the vertex $v_j$ in the deletions of type (T1). Consequently, 
\begin{equation}\label{kanteninE1}
 |E_1|\le \Omega_{i^*+1}\left|\bigcup_{i=i^*+2}^{R+1} X_i(G)\right|k\le \frac{2\Omega_{i^*+1}e(G)}{\Omega_{i^*+2}} \;.
\end{equation}

Consider an arbitrary vertex $w\in\largevertices{\eta}{k}{G}\cap\smallvertices{\eta/2}{k}{G_{j^*}}$ and the interval of those $(j-1)$'s for which $w\in\largevertices{\eta/2}{k}{G_{j-1}} \cap\smallvertices{\eta}{k}{G_{j-1}}$. In such a step the vertex $w$ cannot play the role of the vertices $u_j$ or $v_j$ in~(T2). 
So, each vertex from $\largevertices{\eta}{k}{G}\cap\smallvertices{\eta/2}{k}{G_{j^*}}$ is incident with at least $\eta k/2$ edges from the set $E_0\cup E_1$. Therefore, by the definition of $E_0$, by~\eqref{eq:sparselyconnectedpair}, and by~\eqref{kanteninE1},
$$\left|\largevertices{\eta}{k}{G}\cap\smallvertices{\eta/2}{k}{G_{j^*}}\right|\le \frac{2\cdot|E_0\cup E_1|}{ \eta k/2}\le\left(\frac{4}{R}+\frac{2\Omega_{i^*+1}}{\Omega_{i^*+2}}\right)\cdot\frac{4e(G)}{\eta k}\leByRef{eq:LKSminimalNotManyEdges} \frac{\eta n}2\;.$$
Thus
$$|\largevertices{\eta/2}{k}{G_{j^*}}|\ge |\largevertices{\eta}{k}{G}|-|\largevertices{\eta}{k}{G}\cap\smallvertices{\eta/2}{k}{G_{j^*}}|\ge(1/2+\eta/2)n\;,$$ and consequently, $G_{j^*}\in\LKSgraphs{n}{k}{\eta/2}$. 

Last, we obtain the graph $G'$ by successively deleting any edge from $G_{j^*}$ which connects a vertex from $\smallvertices{\eta/2}{k}{G_{j^*}}$  with a vertex whose degree is not exactly $\lceil(1+\frac\eta2)k\rceil$. This does not affect the already obtained Property~\eqref{item:(prop:gap)gap},
since we could not apply~(T2) to $G_{j^*}$. We claim that for the resulting graph $G'$ we have $G'\in\LKSsmallgraphs{n}{k}{\eta/2}$. Indeed, $\largevertices{\eta/2}{k}{G'}=\largevertices{\eta/2}{k}{G_{j^*}}$, and thus $G'\in\LKSgraphs{n}{k}{\eta/2}$. Property~\ref{def:LKSsmallB} of Definition~\ref{def:LKSsmall} follows from the last step of the construction of $G'$. To see Property~\ref{def:LKSsmallA} of Definition~\ref{def:LKSsmall} we use  Fact~\ref{fact:propertiesOfLKSminimalGraphs}(2)  for $G$ (which by assumption is in $\LKSmingraphs{n}{k}{\eta}$).
\end{proof}


\subsection{Decomposition of graphs with moderate maximum
degree}\label{ssec:class-black} 

First we introduce some useful notions. We start with dense spots which indicate an accumulation of edges in a sparse graph.

\begin{definition}[\bf \index{general}{dense spot}$(m,\gamma)$-dense spot,
\index{general}{nowhere-dense}$(m,\gamma)$-nowhere-dense]\label{def:densespot} An \emph{$(m,\gamma)$-dense spot} in a graph $G$ is a non-empty bipartite sub\-graph  $D=(U,W;F)$ of  $G$ with
$\density(D)>\gamma$ and $\mindeg (D)>m$. We call $G$
\emph{$(m,\gamma)$-nowhere-dense} if it does not contain any $(m,\gamma)$-dense spot.
\end{definition}
We remark that dense spots as bipartite graphs do not have a 
specified orientation, that is, we view $(U,W;F)$ and $(W,U;F)$ as
the same object.

\begin{fact}\label{fact:sizedensespot}
Let $(U,W;F)$ be a $(\gamma k,\gamma)$-dense spot in a
graph $G$ of maximum degree at most $\Omega k$. Then
$\max\{|U|,|W|\}\le \frac{\Omega}{\gamma}k.$
\end{fact}
\begin{proof}
It suffices to observe that $$\gamma |U||W|\leq
e(U,W)\leq\maxdeg(G)\cdot\min\{|U|,|W|\}\leq\Omega k\cdot \min\{|U|,|W|\}.$$
\end{proof}

The next fact asserts that in a bounded degree graph there cannot be too many edge-disjoint dense spots containing a given vertex.
\begin{fact}\label{fact:boundedlymanyspots}
Let $H$ be a graph of maximum degree at most $\Omega k$, let $v\in V(H)$, and let $\DenseSpots$ be a family of edge-disjoint $(\gamma k,\gamma)$-dense spots in $H$. Then less than $\frac{\Omega}{\gamma}$ dense spots from $\DenseSpots$ contain $v$.
\end{fact}
\begin{proof}
This follows as $v$ sends more than $\gamma k$ edges to each dense spot from $\DenseSpots$ it is incident with, the dense spots $\DenseSpots$ are edge-disjoint, and $\deg(v)\le \Omega k$. 
\end{proof}

\medskip
Our second definition of this section might seem less intuitive at first sight.
It describes a property for finding dense spots outside some ``forbidden'' set
$U$, which in later applications will be the set of vertices already used for
a partial embedding of a tree $T\in\treeclass{k}$ from Theorem~\ref{thm:main} during our sequential embedding procedure. In Section~\ref{sssec:whyavoiding} we give a non-technical description of this embedding technique. Informally, a set $\smallatoms$ of vertices is avoiding if for each set $U$ of size $\Theta(k)$ and each vertex $v\in\smallatoms$ there is a dense spot containing $v$ that is almost disjoint from $U$.
\begin{definition}[\bf
\index{general}{avoiding (set)}$(\Lambda,\epsilon,\gamma,k)$-avoiding set]\label{def:avoiding} Suppose that
$G$ is a graph and $\DenseSpots$ is a family of dense spots in $G$. A set
$\smallatoms\subset \bigcup_{D\in\DenseSpots} V(D)$ is \emph{$(\Lambda,\epsilon,\gamma,k)$-avoiding} with
respect to $\DenseSpots$ if for every $U\subset V(G)$ with $|U|\le \Lambda k$ the following holds for all but at most $\epsilon k$ vertices $v\in\smallatoms$. There is a dense spot $D\in\DenseSpots$ with $|U\cap V(D)|\le \gamma^2 k$ that contains $v$.
\end{definition}

Note that a subset of a $(\Lambda,\epsilon,\gamma,k)$-avoiding set is also $(\Lambda,\epsilon,\gamma,k)$-avoiding.

We now come to the main concepts of this section, the bounded and the sparse
decompositions. These notions in a way correspond to the partition structure
from the regularity lemma, although naturally more
complex since we deal with (possibly) sparse graphs here. Lemma~\ref{lem:decompositionIntoBlackandExpanding} is then a corresponding
regularization result. 

\begin{definition}[\index{general}{bounded decomposition}{\bf
$(k,\Lambda,\gamma,\epsilon,\nu,\rho)$-bounded decomposition}]\label{bclassdef}
Suppose that $k\in\NN$ and $\epsilon,\gamma,\nu,\rho>0$ and $\Lambda>2$.
Let $\mathcal V=\{V_1, V_2,\ldots, V_s\}$ be a partition of the vertex set of a graph $G$. We say that $( \clusters,\DenseSpots, \Gblack, \Gexp,
\smallatoms )$ is a {\em $(k,\Lambda,\gamma,\epsilon,\nu,\rho)$-bounded
decomposition} of $G$ with respect to $\mathcal V$ if the following properties
are satisfied:
\begin{enumerate}
\item\label{defBC:nowheredense}
$\Gexp$ is  a $(\gamma k,\gamma)$-nowhere-dense subgraph of $G$ with $\mindeg(\Gexp)>\rho k$.
\item\label{defBC:clusters} The elements of $\clusters$ are disjoint subsets of 
$ V(G)$.
\item\label{defBC:RL} $\Gblack$ is a subgraph of $G-\Gexp$ on the vertex set $\bigcup \clusters$. For each edge
 $xy\in E(\Gblack)$ there are distinct $C_x\ni x$ and $C_y\ni y$ from $\clusters$,
and  $G[C_x,C_y]=\Gblack[C_x,C_y]$. Furthermore, 
$G[C_x,C_y]$ forms an $\epsilon$-regular pair of  density at least $\gamma^2$.
\item\label{bcdef:clustersize} We have $\nu k\le |C|=|C'|\le \epsilon k$ for all
$C,C'\in\clusters$.\label{Csize}
\item\label{defBC:densepairs}  $\DenseSpots$ is a family of edge-disjoint $(\gamma
k,\gamma)$-dense spots  in $G-\Gexp$.  For
each $D=(U,W;F)\in\DenseSpots$ all the edges of $G[U,W]$ are covered
by $\DenseSpots$ (but not necessarily by $D$).
\item\label{defBC:dveapul} If  $\Gblack$
contains at least one edge between $C_1,C_2\in\clusters$ then there exists a dense
spot $D=(U,W;F)\in\DenseSpots$ for which $C_1\subset U$ and $C_2\subset
W$.
\item\label{defBC:prepartition}
For
each $C\in\clusters$ there is a $V\in\mathcal V$ so that either $C\subseteq V\cap V(\Gexp)$ or $C\subseteq V\setminus V(\Gexp)$.
For
each $C\in\clusters$ and $D=(U,W; F)\in\DenseSpots$ we have $C\cap U,C\cap W\in\{\emptyset, C\}$.
\item\label{defBC:avoiding}
$\smallatoms$ is a $(\Lambda,\epsilon,\gamma,k)$-avoiding subset  of
$V(G)\setminus \bigcup \clusters$ with respect to dense spots $\DenseSpots$.
\end{enumerate}

\smallskip
We say that the bounded decomposition $(\clusters,\DenseSpots, \Gblack, \Gexp,
\smallatoms )$ {\em respects the avoiding threshold~$b$}\index{general}{avoiding threshold}\index{general}{respect avoiding threshold} if for each $C\in \clusters$ we either have $\maxdeg_G(C,\smallatoms)\le b$, or $\mindeg_G(C,\smallatoms)> b$.
\end{definition}
Here ``exp'' in $\Gexp$ stands for ``expander'' and ``reg'' in $\Gblack$ stands for ``regular(ity)''.

The members of $\clusters$ are called \index{general}{cluster}{\it clusters}. Define the \index{general}{cluster graph}
{\it cluster graph} \index{mathsymbols}{*Gblack@$\BGblack$}  $\BGblack$ as the graph
on the vertex set $\clusters$ that has an edge $C_1C_2$
for each pair $(C_1,C_2)$ which has density at least $\gamma^2$ in the graph
$\Gblack$. 

Property~\ref{defBC:prepartition} tells us that the clusters may be prepartitioned, just as it is the case in the classic regularity lemma. When in Lemma~\ref{lem:LKSsparseClass} below we classify the graph $G$ from Theorem~\ref{thm:main}   we shall use the prepartition into (roughly) $\smallvertices{\alpha_\PARAMETERPASSING{T}{thm:main}}{k}{G}$ and $\largevertices{\alpha_\PARAMETERPASSING{T}{thm:main}}{k}{G}$.

As said above, the notion of bounded decomposition is needed for our regularity
lemma type decomposition given in
Lemma~\ref{lem:decompositionIntoBlackandExpanding}. It turns out that such a
decomposition is possible only when the graph is of moderate maximum degree. On
the other hand, Lemma~\ref{prop:gapIllustration} tells us that the vertex set of any
graph can be decomposed into vertices of enormous degree and moderate degree.
The graph induced by the latter type of vertices then admits the decomposition from
Lemma~\ref{lem:decompositionIntoBlackandExpanding}. Thus, it makes sense to
enhance the structure of bounded decomposition by vertices of unbounded degree.
This is done in the next definition.

\begin{definition}[\bf \index{general}{sparse
decomposition}$(k,\Omega^{**},\Omega^*,\Lambda,\gamma,\epsilon,\nu,\rho)$-sparse decomposition]\label{sparseclassdef}
Suppose that $k\in\NN$ and $\epsilon,\gamma,\nu,\rho>0$ and $\Lambda,\Omega^*,\Omega^{**}>2$.
Let $\mathcal V=\{V_1, V_2,\ldots, V_s\}$ be a partition of the vertex set of a graph $G$. We say that 
$\class=(\HugeVertices, \clusters,\DenseSpots, \Gblack, \Gexp, \smallatoms )$
is a
{\em $(k,\Omega^{**},\Omega^*,\Lambda,\gamma,\epsilon,\nu,\rho)$-sparse decomposition} of $G$
with respect to $V_1, V_2,\ldots, V_s$ if the following holds.
\begin{enumerate}
\item\label{def:classgap} $\HugeVertices\subset V(G)$,
$\mindeg_G(\HugeVertices)\ge\Omega^{**}k$,
$\maxdeg_H(V(G)\setminus \HugeVertices)\le\Omega^{*}k$, where $H$ is spanned by the edges of $\bigcup\DenseSpots$, $\Gexp$, and
edges incident with $\HugeVertices$,
\item \label{def:spaclahastobeboucla} $( \clusters,\DenseSpots, \Gblack,\Gexp,\smallatoms)$ is a 
$(k,\Lambda,\gamma,\epsilon,\nu,\rho)$-bounded decomposition of
$G-\HugeVertices$ with respect to $V_1\setminus \HugeVertices, V_2\setminus \HugeVertices,\ldots, V_s\setminus \HugeVertices$.
\end{enumerate}
\end{definition}

If the parameters do not matter, we call $\class$ simply a {\em sparse
decomposition}, and similarly we speak about a {\em bounded decomposition}. 

 \begin{definition}[\bf \index{general}{captured edges}captured edges]\label{capturededgesdef}
In the situation of Definition~\ref{sparseclassdef}, we refer to the edges in
$ E(\Gblack)\cup E(\Gexp)\cup
E_G(\HugeVertices,V(G))\cup E_{\GD}(\smallatoms,\smallatoms\cup \bigcup \clusters)$
as \index{general}{captured edges}{\em captured} by the sparse decomposition. 
 We write
\index{mathsymbols}{*Gclass@$\Gcapt$}$\Gcapt$ for the subgraph of $G$ on the same vertex set which consists of the captured edges.
Likewise, the captured edges of a bounded decomposition
$(\clusters,\DenseSpots, \Gblack,\Gexp,\smallatoms )$ of a graph $G$ are those
in $E(\Gblack)\cup E(\Gexp)\cup E_{\GD}(\smallatoms,\smallatoms\cup\bigcup\clusters)$.
\end{definition}

Throughout the paper we write \index{mathsymbols}{*GD@$\GD$}$\GD$ for the subgraph of $G$
which consists of the edges contained in $\DenseSpots$. We now include an easy fact about the relation of $\GD$ and $\Gblack$.
\begin{fact}\label{fact:denseVSblack}
Let $\class=(\HugeVertices, \clusters,\DenseSpots, \Gblack, \Gexp,\smallatoms )$ be a sparse decomposition of
a graph $G$. Then each edge $xy\in E(\GD)$ with
$x,y\in\bigcup\clusters$ is either contained in $\Gblack$, or is not captured.
\end{fact}
\begin{proof}
Indeed, suppose that $xy\in E(\GD)$,
$x,y\in\bigcup\clusters$, and $xy\not\in E(\Gblack)$.  Property~\ref{def:spaclahastobeboucla} of
Definition~\ref{sparseclassdef} says that $x,y\notin \HugeVertices$. Further, by Property~\ref{defBC:avoiding} of Definition~\ref{bclassdef}, we
have $x,y\not\in \smallatoms$.  Last, Property~\ref{defBC:densepairs} of Definition~\ref{bclassdef} implies that $xy\not\in E(\Gexp)$. Hence $xy$ is not captured, as desired.
\end{proof}

We now give a bound on the number of clusters reachable through edges of
the dense spots from a fixed vertex outside $\HugeVertices$. 

\begin{fact}\label{fact:clustersSeenByAvertex}
Let  $\class=(\HugeVertices, \clusters,\DenseSpots, \Gblack, \Gexp,\smallatoms )$ be a 
$(k,\Omega^{**},\Omega^*,\Lambda,\gamma,\epsilon,\nu,\rho)$-sparse
decomposition of a graph  $G$. Let $x\in V(G)\setminus \HugeVertices$. Assume that $\clusters\not=\emptyset$, and let $\clustersize$ be the size of each of the members of~$\clusters$. Then there are less than
$$\frac{2(\Omega^*)^2k}{\gamma^2 \clustersize}\le\frac{2(\Omega^*)^2}{\gamma^2\nu}$$ clusters
$C\in\clusters$ with $\deg_{\GD}(x,C)>0$.
\end{fact}
\begin{proof}
Property~\ref{def:classgap} of Definition~\ref{sparseclassdef} says that
$\deg_{\GD}(x)\le \Omega^*k$. For each $D\in
\DenseSpots$ with $x\in V(D)$
we have that
$\deg_{D}(x)>
\gamma k$, since $D$ is a $(\gamma k,\gamma )$-dense spot. 
By Fact~\ref{fact:boundedlymanyspots}
\begin{equation}\label{labello}
 |\{D\in\DenseSpots:\deg_{D}(x)>0\}|< \frac {\Omega^*}{\gamma}.
\end{equation}

Furthermore, 
by Fact~\ref{fact:sizedensespot}, and using Properties~\ref{Csize} and~\ref{defBC:dveapul} of Definition~\ref{bclassdef}, we see that for a fixed~$D\in\DenseSpots$, we have
$$|\{C\in\clusters\::\:
C\subset V(D)\}|\le \frac{2\Omega^*
k}{\gamma}\cdot\frac1{\clustersize}\le\frac{2\Omega^*}{\gamma\nu}\;.$$
Together with~\eqref{labello} this gives that the number of clusters $C\in\clusters$ with
$\deg_{\GD}(x,C)>0$ is less than $$\frac{\Omega^*}{\gamma}\cdot
\frac{2\Omega^*k}{\gamma \clustersize}\le\frac{\Omega^*}{\gamma}\cdot
\frac{2\Omega^*}{\gamma\nu}\;,$$ as desired.
\end{proof}

As a last step before we state the main result of this section we show that the
cluster graph $\BGblack$ corresponding to a
$(k,\Omega^{**},\Omega^*,\Lambda,\gamma,\epsilon,\nu,\rho)$-sparse
decomposition $(\HugeVertices, \clusters,\DenseSpots, \Gblack, \Gexp,\smallatoms)$  has bounded degree.

\begin{fact}\label{fact:ClusterGraphBoundedDegree}
Let  $\class=(\HugeVertices, \clusters,\DenseSpots, \Gblack, \Gexp,\smallatoms )$ be a 
$(k,\Omega^{**},\Omega^*,\Lambda,\gamma,\epsilon,\nu,\rho)$-sparse
decomposition of a graph  $G$, and let $\BGblack$ be the
corresponding cluster graph. Let $\clustersize$ be the size of each cluster in $\clusters$. 
Then $\maxdeg(\BGblack)\le \frac{\Omega^*
k}{\gamma^2\clustersize}\le\frac{\Omega^*}{\gamma^2\nu}$.
\end{fact}
\begin{proof}
Let $C\in\clusters$. Then by the definition of $\BGblack$ and by Property~\ref{defBC:RL} of
Definition~\ref{bclassdef} we have $\deg_{\BGblack}(C)\leq \sum_{C'\in\neighbour_{\BGblack}(C)}\frac{e_{\Gblack}(C,C')}{\gamma^2|C||C'|}=\sum_{C'\in\neighbour_{\BGblack}(C)}\frac{e_{\Gblack}(C,C')}{\gamma^2|C|\clustersize}$. Since the maximum degree in $\Gblack$ is upper-bounded by $\Omega^*k$ (c.f.\ Property~\ref{def:classgap} of Definition~\ref{sparseclassdef}), we get 
$$\deg_{\BGblack}(C)\leq \sum_{C'\in\neighbour_{\BGblack}(C)}\frac{e_{\Gblack}(C,C')}{\gamma^2|C|\clustersize} \leq\frac{\Omega^*k|C|}{ \gamma^2 |C|\clustersize}\; \leBy{\mbox{D\ref{bclassdef}(\ref{bcdef:clustersize})}}\; \frac{\Omega^*}{ \gamma^2 \nu }\;,$$
 as desired.
\end{proof}

\medskip

We now state the most important lemma of this section. It says that any
 graph of bounded degree has a bounded decomposition which captures almost all its edges. This lemma can be considered as a sort
of regularity lemma for sparse graphs.

\begin{lemma}[Decomposition lemma]\label{lem:decompositionIntoBlackandExpanding}
For each $\Lambda,\Omega,s\in\NN$  and each $\gamma,\epsilon,\rho>0$ there
exist $k_0\in\NN$, $\nu>0$ such that for every $k\ge k_0$ and every
$n$-vertex graph $G$ with $e(G)\le kn$, $\maxdeg(G)\le \Omega k$, and with a
given partition $\mathcal V$ of its vertex set into at most $s$ sets, the following holds for each $b>0$. There
exists a $(k,\Lambda,\gamma,\epsilon,\nu,\rho)$-bounded decomposition
$(\clusters,\DenseSpots, \Gblack, \Gexp,\smallatoms )$ with respect to $\mathcal
V$, which captures all but at most
$(\frac{4\epsilon}{\gamma}+\epsilon\Omega+\gamma+\rho)kn$ edges of $G$ and respects avoiding threshold $b$. 
Furthermore, we have
\begin{equation}\label{eq:sEr}
 |E(\DenseSpots)\setminus (E(\Gblack)\cup
E_{\GD}[\smallatoms,\smallatoms\cup\bigcup\clusters])|\le
(\frac{4\epsilon}{\gamma}+\epsilon\Omega +\gamma) kn\;.
\end{equation}
\end{lemma}

A proof of Lemma~\ref{lem:decompositionIntoBlackandExpanding} is given in
Section~\ref{ssec:ProofOfDecomposition}. 

\subsection{Decomposition of LKS graphs}\label{ssec:DecompositionOfLKSGraphs}
Lemma~\ref{prop:gap} and
Lemma~\ref{lem:decompositionIntoBlackandExpanding} enable us to decompose graphs
in $\LKSgraphs{n}{k}{\eta}$ in a particular manner.

\begin{lemma}\label{lem:LKSsparseClass}
For every $\eta, \Lambda,\gamma,\epsilon,\rho\in(0,1)$ there are $\nu>0$ and $k_0\in\NN$ such
that for every $k>k_0$ and  for every number $b>0$ the following holds.
For every   sequence $(\Omega_j)_{j\in\NN}$ of positive numbers with
$\Omega_1>2$, $\Omega_j/\Omega_{j+1}\le \eta^2/100$ for all $j\in\NN$
 and for every $G\in\LKSgraphs{n}{k}{\eta}$ there are an index $i$ and a subgraph $G' $ of $G$ with the following properties:
 \begin{enumerate}[(a)]
 \item $G'\in\LKSsmallgraphs{n}{k}{\eta/2}$,
 \item $i\leq 100\eta^{-2}$,
 \item\label{LKSclassif:prepart} $G'$ has a
 $(k,\Omega_{i+1},\Omega_{i},\Lambda,\gamma,\epsilon,\nu,\rho)$-sparse decomposition
$(\HugeVertices, \clusters,\DenseSpots,
\Gblack',\Gexp',\smallatoms)$ with respect to the partition
$\{V_1,V_2\}:=\{\smallvertices{\eta/2}{k}{G'},\largevertices{\eta/2}{k}{G'}\}$, and with respect to avoiding threshold $b$,
\item\label{it:LKSsparseGreyCapt} $(\HugeVertices, \clusters,\DenseSpots,
\Gblack',\Gexp',\smallatoms)$ captures all but at most
$(\frac{4\epsilon}{\gamma}+\epsilon\Omega_{\lfloor
100\eta^{-2}\rfloor}+\gamma+\rho )kn$ edges of $G'$, and
\item $|E(\DenseSpots)\setminus (E(\Gblack')\cup
E_{G'}[\smallatoms,\smallatoms\cup\bigcup\clusters])|\le
(\frac{4\epsilon}{\gamma}+\epsilon\Omega_{\lfloor
100\eta^{-2}\rfloor}+\gamma) kn$.
\end{enumerate}
\end{lemma}
\begin{proof}
Let $\nu$ and $k_0$
be given by Lemma~\ref{lem:decompositionIntoBlackandExpanding}
for input parameters
$\Omega_\PARAMETERPASSING{L}{lem:decompositionIntoBlackandExpanding}:=\Omega_{\lfloor 100\eta^{-2}\rfloor}$,
$\Lambda_\PARAMETERPASSING{L}{lem:decompositionIntoBlackandExpanding}:=\Lambda,
\gamma_\PARAMETERPASSING{L}{lem:decompositionIntoBlackandExpanding}:=\gamma,
\epsilon_\PARAMETERPASSING{L}{lem:decompositionIntoBlackandExpanding}:=\epsilon,
\rho_\PARAMETERPASSING{L}{lem:decompositionIntoBlackandExpanding}:=\rho,
b_\PARAMETERPASSING{L}{lem:decompositionIntoBlackandExpanding}:=b$, 
 and
$s_\PARAMETERPASSING{L}{lem:decompositionIntoBlackandExpanding}:=2$. Now, given
$G$, let us consider a subgraph $\tilde G$ of $G$ such that $\tilde G\in
\LKSmingraphs{n}{k}{\eta}$. Lemma~\ref{prop:gap} applied to the sequence
$(\Omega_j)_j$ and $\tilde G$ yields a graph $G'\in \LKSsmallgraphs{n}{k}{\eta/2}$
and an index $i\leq 100\eta^{-2}$. We set $\HugeVertices:=\{v\in V(G)\::\:\deg_{G'}(v)\ge
\Omega_{i+1}k\}$.

Observe that by~\eqref{eq:LKSminimalNotManyEdges}, $e(G')<kn$. Let
$(\HugeVertices,\DenseSpots, \Gblack',\Gexp',\smallatoms)$ be the
$(k,\Lambda,\gamma,\epsilon,\nu,\rho)$-bounded decomposition of the graph
$G'-\HugeVertices$ with respect to $\{\smallvertices{\eta/2}{k}{G'},\largevertices{\eta/2}{k}{G'}\setminus \HugeVertices\}$ that
is given by Lemma~\ref{lem:decompositionIntoBlackandExpanding}. Clearly,
$(\HugeVertices, \clusters,\DenseSpots,
\Gblack',\Gexp',\smallatoms)$ is  a
$(k,\Omega_{i+1},\Omega_{i},\Lambda,\gamma,\epsilon,\nu,\rho)$-sparse
decomposition of~$G'$ capturing at least as many edges as promised in the lemma.
\end{proof}

\bigskip
The process of embedding a given tree $T_\PARAMETERPASSING{T}{thm:main}\in\treeclass{k}$ into $G_\PARAMETERPASSING{T}{thm:main}$ is based on the sparse decomposition
$\class=(\HugeVertices, \clusters,\DenseSpots,\Gblack,\Gexp,\smallatoms)$ of~$G_\PARAMETERPASSING{T}{thm:main}$ given by Lemma~\ref{lem:LKSsparseClass} and is much more complex than in
approaches based on the standard regularity lemma. The embedding ingredient in the classic (dense)
regularity method inheres in blow-up lemma type statements which roughly tell
that  regular pairs of positive density in some sense behave like complete bipartite graphs.
In our setting, in addition to regular pairs
 we shall use three other components of $\class$: the vertices of
huge degree $\HugeVertices$, the nowhere-dense graph $\Gexp$, and the avoiding
set $\smallatoms$. Each of these components requires a different strategy for
embedding (parts of) $T_\PARAMETERPASSING{T}{thm:main}$. Let us mention that rather major technicalities arise when combining these strategies.

These strategies are described precisely and in detail
in~\cite{cite:LKS-cut3}. An informal account on the role of
$\smallatoms$ is given in Section~\ref{sssec:whyavoiding}. We discuss the use of $\Gexp$ in
Section~\ref{sssec:whyGexp}. Only very little can be said about the
set $\HugeVertices$ at an intuitive level: these vertices have huge degrees but
are very unstructured otherwise. If  only $o(kn)$ edges  are incident
with $\HugeVertices$ then we can neglect them. If, on the other hand, there are  $\Omega(kn)$ edges incident with 
$\HugeVertices$, then we have no choice but to use them for our embedding. Very roughly speaking, in that case we
find sets $\HugeVertices'\subset\HugeVertices$ and $V'\subset V(G)\setminus
\HugeVertices$ such that still $\mindeg(\HugeVertices',V')\gg k$, and
$\mindeg(V',\HugeVertices')=\Omega(k)$, and then use $\HugeVertices'$ and $V'$ in our
embedding.

\medskip
Last, let us note that when $G_\PARAMETERPASSING{T}{thm:main}$ is close to the extremal graph (depicted in Figure~\ref{fig:ExtremalGraph}) then all the structure in $G_\PARAMETERPASSING{T}{thm:main}$ captured by Lemma~\ref{lem:LKSsparseClass} accumulates in the cluster graph $\Gblack'$, i.e., $\HugeVertices$, $\Gexp'$ and $\smallatoms$ are all almost empty. For that reason, when some of $\HugeVertices$, $\Gexp'$ or $\smallatoms$ is substantial we gain some extra aid. In comparison, one of the almost extremal graphs for the Erd\H{o}s--S\'os Conjecture~\ref{conj:ES} has a substantial $\HugeVertices$-component (see Figure~\ref{fig:ExtremalGraphES2}).

\subsection{Decomposition of general graphs}
A version of Lemma~\ref{lem:LKSsparseClass} can be formulated for general graphs. To illustrate this, we present below a generic lemma of this type, which will not be used in the proof of the main theorem. 

\begin{lemma}\label{lem:genericBD}
For every $\eta, \Lambda,\gamma,\epsilon,\rho>0$ there are numbers $\nu>0$ and $k_0\in\NN$ such that
for every   sequence $(\Omega_j)_{j\in\NN}$ of positive numbers with
$\frac{\Omega_j}{\Omega_{j+1}}\le \frac{\eta}4$ the following holds. Suppose that $G$ is a graph of order $n$ with average degree $k>k_0$. Then there is an index $i\leq \frac4\eta$, such that 
$G$ has a
 $(k,\Omega_{i+1},\Omega_{i},\Lambda,\gamma,\epsilon,\nu,\rho)$-sparse decomposition
$(\HugeVertices, \clusters,\DenseSpots,
\Gblack,\Gexp,\smallatoms)$ 
that captures all but at most
\begin{equation}
\label{eq:generalDecUncapt}
(\eta+\frac{4\epsilon}{\gamma}+\epsilon\Omega_{\lfloor
4\eta^{-1}\rfloor}+\gamma+\rho )kn
\end{equation}
edges.
\end{lemma}
The proof follows the same strategy as that of Lemma~\ref{lem:LKSsparseClass}.
\begin{proof}[Proof outline]
By Lemma~\ref{prop:gapIllustration} there exists a spanning sugraph $G'$ of $G$ with $e(G)-e(G')<\eta kn$, and an index $i\leq \frac4\eta$ such that the assertion of Lemma~\ref{prop:gap}\eqref{item:(prop:gap)gap} holds. The bounded-degree part can then be decomposed using Lemma~\ref{lem:decompositionIntoBlackandExpanding}, yielding the desired sparse decomposition.
\end{proof}
This decomposition could be used to attack other problems; probably with a version of Lemma~\ref{lem:genericBD} tailored to a particular setting similarly as we did in Lemma~\ref{lem:LKSsparseClass}.
 However, our feeling is that such a decomposition lemma is limited in applications to tree-containment problems. The reason is that two of the features of the sparse decomposition, the nowhere-dense graph $\Gexp$ and the avoiding set $\smallatoms$, seem to be useful only for embedding trees. See Section~\ref{sssec:whyavoiding} and Section~\ref{sssec:whyGexp} for a discussion of the respective embedding strategies.

\subsection{The role of the avoiding set $\smallatoms$}\label{sssec:whyavoiding}
Let us explain the role of the avoiding set $\smallatoms$
in Lemma~\ref{lem:decompositionIntoBlackandExpanding}. As said above, our aim in
Lemma~\ref{lem:decompositionIntoBlackandExpanding} will be to locally regularize
parts of the input graph $G$. Of course, first we try to regularize as large a part of the
$G$ as possible. The avoiding set arises as a result of the
impossibility to regularize certain parts of the graph. Indeed, it is one of the most surprising steps in our proof of
Theorem~\ref{thm:main} that the set $\smallatoms$ is initially defined as ---
very loosely speaking --- ``those vertices where the regularity lemma fails to
work properly'', and only then we prove
that $\smallatoms$
actually satisfies the useful conditions of Definition~\ref{def:avoiding}.

We now sketch how to utilize avoiding sets for the purpose of embedding trees.
 In our proof of Theorem~\ref{thm:main} we preprocess the tree
$T=T_\PARAMETERPASSING{T}{thm:main}\in\treeclass{k}$ by choosing several cut-vertices so that the tree decomposes into small components, called \emph{shrubs}. We cut $T$ so that the order of each shrub is at most $\tau k$, where $\tau>0$ is a small constant. 
Then we
sequentially embed those shrubs. Thus embedding techniques for
embedding a single shrub are the building blocks of our embedding machinery; and $\smallatoms$ is one of the enviroments which
provides us with such a technique. Let us discuss here the simpler case of embedding end
shrubs (i.e.\ shrubs incident to a single cut-vertex). More precisely, we show how to extend a partial embedding of a tree by one end-shrub. To this end, let us suppose that $\phi$ is a partial embedding of a tree $T$, and $v\in V(T)$ is its \emph{active vertex}\index{general}{active vertex}, i.e., a vertex which is embedded, but not all its children are. We write $ U\subset V(G)$ for the current image of $\phi$. Let
$T'\subset T$ be an end-shrub which is not embedded yet, and suppose $u\in V(T')$ is adjacent to $v$. We have $v(T')\le\tau k$. 

\begin{figure}[ht] \centering
\includegraphics[scale=1.00]{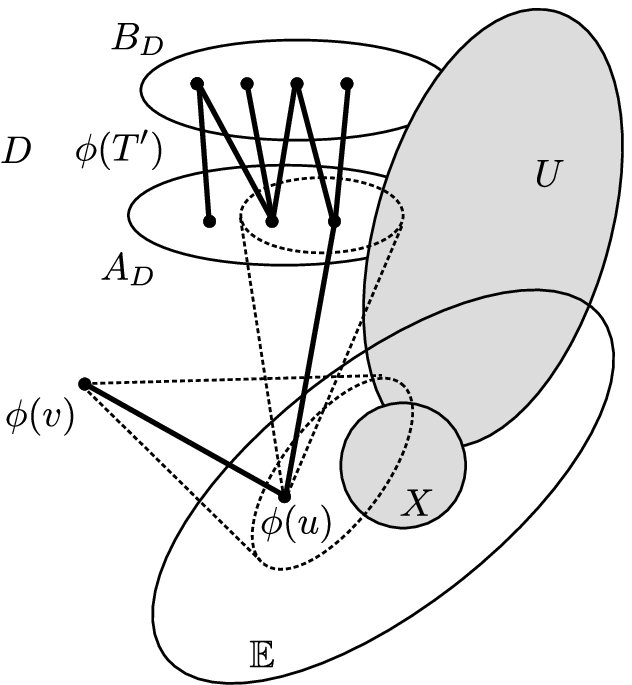}
\caption[Embedding using the set $\smallatoms$]{Embedding using the set $\smallatoms$.}
\label{fig:EmbeddingAvoiding}
\end{figure}

We now show how to
extend the partial embedding $\phi$ to $T'$, assuming that
$\deg_G\big(\phi(v),\smallatoms\setminus U\big)\ge \gamma k$ for some
$(1,\epsilon,\gamma,k)$-avoiding set $\smallatoms$ (where $\tau\ll
\epsilon\ll \gamma\ll 1$). Let $X$ be the set of at most $\epsilon k$
exceptional vertices from Definition~\ref{def:avoiding} corresponding to the set
$U$. We now embed $T'$ into $G$, starting by embedding $u$ in a vertex of
$\smallatoms\setminus (U\cup X)$ in the neighbourhood of $\phi(v)$. By Definition~\ref{def:avoiding}, there is a dense spot
$D=(A_D,B_D; F)\in\DenseSpots$ such that $\phi(u)\in V(D)$
and $|U\cap V(D)|\le \gamma^2k$. 
As $D$ is a dense spot, we have $\deg_G(\phi(u),V(D))>\gamma k$.
We can greedily embed $T'$ into $D$ using the minimum degree in $D$. See
Figure~\ref{fig:EmbeddingAvoiding} for an illustration, and \cite[Lemma~\ref{p3.lem:embed:avoidingFOREST}]{cite:LKS-cut3}
for a
precise formulation.

We indeed use the avoiding set for embedding shrubs of
$T$ as above. The major simplification we made in the exposition is that we
only discussed the case when $T'$ is an end shrub. To cover embedding of an internal
shrub $T'$ as well (i.e.\ a shrub that is incident to more than one cut-vertex), one needs to have a more detailed control over the embedding, i.e., one must be able to extend the embedding of $T'$ to the neighbouring cut-vertices, 
in such a way that one can then continue the embedding. 

Last, let us remark, that unlike our baby-example above, we use an
$(\Lambda,\epsilon,\gamma,k)$-avoiding set with $\Lambda\gg 1$. This is because
in the actual proof one has to avoid more vertices than just the current image
of the embedding.

\subsection{The role of the nowhere-dense graph $\Gexp$}\label{sssec:whyGexp}

In this section we shall give some intuition on how the $(\gamma k, \gamma)$-nowhere-dense graph $\Gexp$ from the
$(k,\Omega^{**},\Omega^*,\Lambda,\gamma,\epsilon',\nu,\rho)$-sparse
decomposition\footnote{We shall assume that $17\sqrt{\gamma}<\rho$; this will be
the setting of the sparse decomposition we shall work with in the proof of
Theorem~\ref{thm:main}.} $(\HugeVertices, \clusters,\DenseSpots, \Gblack,
\Gexp,\smallatoms)$ of a graph $G$ is useful for embedding a
given tree $T\in\treeclass{k}$. We start out with the rather simple case when $T$ is a path. We then  point out an issue
with this approach for trees with many branching vertices and show how to overcome this problem. 

\paragraph{Embedding a path in $\Gexp$.} 
Assume we are given a path $T=u_1u_2\cdots u_k\in\treeclass{k}$ and we wish to embed it into $\Gexp$.
The idea is to apply a one-step look-ahead strategy. We first embed  $u_1$ in an arbitrary vertex $v\in V(\Gexp)$. Then, we extend our
embedding $\phi_\ell$ of the path $u_1\cdots u_\ell$ in $\Gexp$ in 
step $\ell$
 by em\-bedding $u_{\ell +1}$ in a (yet unused) neighbour $w$ of the image of the 
 \emph{active} vertex~$u_\ell$, requiring that
\begin{equation}\label{eq:inductiveembedding3}
\deg_{\Gexp}\big(w,\phi_\ell(u_1\cdots u_\ell)\big)<\sqrt\gamma k\;.
\end{equation}
Let us argue that such a vertex $w$ exists using induction on $\ell$. First, observe that
Property~\ref{defBC:nowheredense} of Definition~\ref{bclassdef} implies that  $\phi_\ell(u_\ell)$ has  at least $\rho k$ neighbours.
By~\eqref{eq:inductiveembedding3} applied to $\ell -1$, at most $\sqrt\gamma k$ of these neighbours
lie inside $\phi_\ell(u_1\cdots u_{\ell-1})$; this property is also trivially satisfied when $\ell=1$. 
Further, an easy calculation shows that at most $16\sqrt\gamma k$ of them have degree more than $\sqrt\gamma k$ in $\Gexp$ into the set $\phi_\ell(u_1\cdots u_\ell)$, otherwise we would get a contradiction to $\Gexp$ being $(\gamma k, \gamma)$-nowhere-dense. Since we assumed $\rho> 17\sqrt\gamma$ we can find a vertex $w$ satisfying~\eqref{eq:inductiveembedding3} and thus embed all of~$T$.


\paragraph{Embedding trees with many branching points.}
We certainly cannot hope that a nonempty graph $\Gexp$ alone will provide us
with embeddings of all trees $T\in\treeclass{k}$ from Theorem~\ref{thm:main}. For instance, if $T$ is a star, then we need in $G$ a vertex of degree $k-1$, which $\Gexp$ might not have. The structure of LKS graphs allows to deal with embedding high-degree vertices. However, even without any vertex of large degree in our tree, the method described above might not always work, as we show next.

\begin{figure}[t]
\centering 
\includegraphics{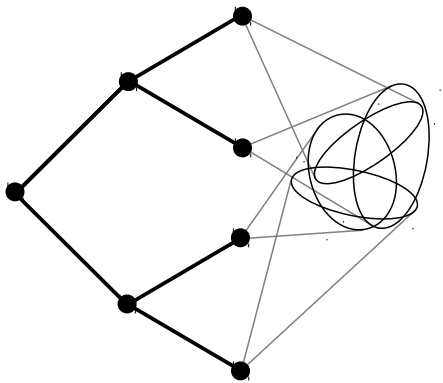}
\caption[Getting stuck while embedding binary tree]{Embedded part of the binary tree in bold. 
The neighbourhoods of active vertices may overlap.}\label{fig:gettingstuck}
\end{figure}

Consider a binary tree $T\in\treeclass{k}$, rooted
at its central vertex $r$. Now if we try to embed~$T$
 sequentially as above we will arrive at a moment
when there are many (as many as $\log_2 k$) active vertices; regardless in which order we embed.\footnote{The only requirement on the ordering is that in each moment the embedded part of the tree forms a connected subgraph; in particular we may use the depth-first and the breadth-first orders.} Now, the
neighbourhoods of the images of the active vertices cannot be controlled much, i.e., they may be  intersecting considerably. Hence, when embedding  children of active vertices we might block available space in the neighbourhoods of other active
vertices. 
 See
Figure~\ref{fig:gettingstuck} for an illustration.

To rescue the situation we partition $T$ so that the first $q$ levels of $T$ from the root $r$
form the set of the cut-vertices $W$. All other vertices
make up the end shrubs $T^*_1,\ldots,T^*_h$. That is, $|W|=2^q-1$, and $h=2^{q+1}-2$.

We first embed the few cut-vertices $W$. As $\rho k$ will be much larger $2^q$, following a strategy similar to the one above we ensure that all cut-vertices get correctly embedded.
The next step is to make the transitions
 at the
 $q$-th level from embedding cut-vertices to embedding shrubs $T^*_1,\ldots,T^*_h$. But since this step requires to exploit the structure of LKS graphs, we skip the details in the high-level overview here. For the sake of this simplified example, let us assume that all the cut-vertices are embedded in a set $L$ with 
 \begin{equation}\label{eq:highdegcut}
\mindeg_{\Gexp}(L)\ge \delta k\; ,
 \end{equation}
 (where $\rho\ll\delta<1$ is a small constant).
 
For the point we wish to make here, it is more relevant to see how to complete the last part of our embedding, that is, how to  embed a tree $T^*_i$ whose root $r_i$
 is already embedded in a vertex $\phi(r_i)\in  V(\Gexp)$.
Let $\text{im}_i:=\text{im}(\phi)$ be the current (partial) image of $\phi$. Further, we assume that throughout the entire process we have
\begin{equation}\label{eq:magic}
\deg_{\Gexp}(\phi(r_i),V(\Gexp)\setminus \text{im})\gtrapprox \delta k/2\;,
\end{equation}
where $\text{im}$ is the image at that moment (and in particular, also at the end of the process).
We explain how to achieve this property at the end. 

We emphasize that at this moment we are working exclusively with the tree $T^*_i$, i.e., any other tree
$T^*_j$ is either completely embedded, or will be embedded only after we finish
the embedding of $T^*_i$. Suppose we are about to embed a vertex $v\in V(T^*_i)$ whose ancestor $v'\in
V(T^*_i)\cup W$ is already embedded in $V(\Gexp)$. We choose for the image of $v$ any (yet unused) vertex $w$ in the neighbourhood of $\varphi (v')$, requiring that 
\begin{equation}\label{eq:degexpi}
\deg_{\Gexp}(w,\text{im}_i)<\rho k/100\;.
\end{equation}
This condition is very similar to our path-embedding procedure above, and can be proved in exactly the same way, using the fact that $\Gexp$ is $(\gamma k, \gamma)$-nowhere-dense. When $v'\in W$ is a cut-vertex, we need to combine this argument with~\eqref{eq:magic}.

Note that during our embedding
$|\text{im}(\phi)\setminus \text{im}_i |$ will grow. However, $|\text{im}(\phi)\setminus \text{im}_i |$ is at most $v(T^*_i)$, which is much smaller than $\rho k$. Thus, for every vertex $v''\in
V(T^*_i)$, when its time to be embedded comes, we still have a small degree 
to the partial image of the tree. Therefore $v''$ can be embedded on a vertex $w$ that satisfies $\deg_{\Gexp}(w,\text{im}(\phi))<\rho k/50$ similarly as in~\eqref{eq:degexpi}.

Note that the trick here was to keep on working on one subtree $T^*_i$, whose size is small enough to be negligible in comparison to the degree of the vertices in $\Gexp$. So, by avoiding the vertices that have a considerable degree into $\text{im}_i$, we actually also avoid those vertices that have a  considerable degree into $\text{im}(\phi)$. 
Breaking up the tree into tiny shrubs was thus the key to successfully embedding it.

\medskip

Let us now explain how to achieve~\eqref{eq:magic}. Instead of just embedding the tree $T^*_i$ we shall also reserve an equal amount of vertices in $\Gexp$ that are touched only exceptionally. More precisely, in a given step, instead of extending the embedding from a vertex to its two children, we first find four candidate vertices, and we randomly select two of them to host these children and insert the remaining two into a reserve set $R$. Condition~\eqref{eq:degexpi} is replaced by $\deg_{\Gexp}(w,\text{im}_i\cup R)<\rho k/100$. This allows us to 
avoid not only $\text{im}_i$ but also $R$ when extending the embedding of $T_i^*$. The only time the set $R$ may be used to host a vertex $v$ of some tree $T_1^*,\ldots,T_h^*$ is when $v$ is the root of such a tree. Since the choice for the inclusion of vertices to $R$ was random, with high probability we have
$$
\deg_{\Gexp}(\phi(r_i),\text{im}_i)\approx \deg_{\Gexp}(\phi(r_i),R)
\pm h\;,$$
where the $\pm h$ term amounts to the roots for which the random choice is not used. Recall that $h\ll \rho k$. This together with~\eqref{eq:highdegcut} establishes~\eqref{eq:magic}.

\subsection{Proof of the decomposition lemma}
\label{ssec:ProofOfDecomposition}
This subsection is devoted to the proof of
the decomposition lemma (Lemma~\ref{lem:decompositionIntoBlackandExpanding}). In the proof, we start by extracting the edges of as many $(\gamma
k,k)$-dense spots from $G$ as possible; these together with the incident vertices will form the auxiliary graph $\GD$. Most of the
remaining edges will form the edge set of the graph $\Gexp$. Next, we consider
the intersections of the dense spots captured in $\GD$. We apply the regularity lemma for locally dense graphs
(Lemma~\ref{lem:sparseRL}) to the subgraph of $\GD$
that is spanned by the large intersections, and thus obtain $\Gblack$. The other part of
$V(\GD)$ will be taken as the  $(\Lambda,\epsilon,\gamma,k)$-avoiding
set~$\smallatoms$.

\smallskip

\paragraph{Setting up the parameters.} We start by setting 
$$\tilde\nu:=\epsilon\cdot 3^{-\frac{\Omega\Lambda}{\gamma^3}}.$$
Let $q_\mathrm{MAXCL}$ be  given by Lemma~\ref{lem:sparseRL} for input parameters 
\begin{equation}\label{eq:in213}
m_\PARAMETERPASSING{L}{lem:sparseRL}:=\frac{\Omega}{\gamma\tilde\nu}\quad, \quad z_\PARAMETERPASSING{L}{lem:sparseRL}:=4s\quad \mbox{ and } \quad
\epsilon_\PARAMETERPASSING{L}{lem:sparseRL}:=\epsilon\;.
\end{equation}
 Define
an auxiliary parameter
$q:=\max\{q_\mathrm{MAXCL}, \epsilon^{-1}\}$ and choose the output parameters of
Lemma~\ref{lem:decompositionIntoBlackandExpanding} as $$k_0:=\left\lceil\frac{
q_\mathrm{MAXCL}}{\tilde\nu}\right\rceil \ \    \ \ \text{ and } \ \    \ \
\nu:=\frac{\tilde \nu}q .$$

\paragraph{Defining $\DenseSpots$ and $\Gexp$.}
Given a graph $G$, take a family $\DenseSpots$ of edge-disjoint
$(\gamma k, \gamma)$-dense spots such that the resulting
graph $\GD\subset G$ (which contains those vertices and edges that are
contained in~$\bigcup\DenseSpots$) has the maximum number  of edges. 

Then by
Lemma~\ref{lem:subgraphswithlargeminimumdegree} there exists a graph
$\Gexp\subset G- \GD$ with $\mindeg(\Gexp)> \rho k$ and such that
\begin{equation}\label{eq:almostalldashed}
|E(G)\setminus (E(\Gexp)\cup  E(\GD))|\le\rho kn\;.
\end{equation} 
This choice of $\DenseSpots$ and $\Gexp$ already satisfies Properties~\ref{defBC:densepairs} and~\ref{defBC:nowheredense}  of Definition~\ref{bclassdef}.


\paragraph{Preparing for an application of the regularity lemma.}
Let
\begin{equation}\label{eq:Xref}
 \mathcal X:= \bigboxplus_D \{ U, W, V(G)\setminus V(D)\}\;,
\end{equation}
where the latter partition refinement ranges over all $D=(U,W;F)\in\DenseSpots$.
Let $\mathcal B:=\{X\in\mathcal X\::\: X\subset V(\GD)\}$, $\tilde{\mathcal B}:=\{B\in\mathcal{B}\::\:|B|>2\tilde\nu k \}$, and 
$\tilde{\mathcal C}:=\mathcal B\setminus\tilde{\mathcal B}$. Furthermore let $\tilde B:=\bigcup_{B\in\tilde{\mathcal B}}B$ and  $\smallatoms:=\bigcup_{C\in\tilde{\mathcal C}}C$.  

Now,  partition each set $B\in\tilde{ \mathcal B}$  into
$c_B:=\lceil|B|/2\tilde\nu k\rceil$ subsets $B_1,\ldots,B_{c_B}$ of
cardinalities differing by at most one, and let $\mathcal B'$ be the set
containing all the sets  $B_i$ (for all $B\in
\tilde{\mathcal B}$).
Then for each $B\in\mathcal B'$ we have that
\begin{equation}\label{eq:ClustersOfRightSize}
\tilde\nu k\le |B| \le 2 \tilde\nu k\le \epsilon k\;.
\end{equation}

Construct a graph $H$ on 
$\mathcal B'$ by making two vertices $A_1,A_2\in \mathcal
B'$ adjacent in $H$ if 
\begin{enumerate}[(A)]
\item\label{it:VV} there is a dense spot $D=(U,W; F)\in \DenseSpots$ such that  $A_1\subset U$ and $A_2\subset W$, and
\item\label{item:density} $\density_{G}(A_1,A_2)\ge \gamma$.
\end{enumerate}
Note that it follows from the way $\DenseSpots$ was chosen that if $A_1A_2\in E(H)$ then $G[A_1,A_2]=\GD[A_1,A_2]$. On the other hand note that we do not necessarily have $G[A_1,A_2]=D[A_1,A_2]$ for the dense spot $D$ appearing in~\eqref{it:VV}; just because there may be several such dense spots $D$.

By the assumption of Lemma~\ref{lem:decompositionIntoBlackandExpanding}, $\maxdeg(G)\le\Omega k$. So, for each
$B\in\mathcal B'$ we have
$e_G(B,\tilde B\setminus B)\le
\Omega k|B|$. On the other hand,
\eqref{eq:ClustersOfRightSize} and~\eqref{item:density}
imply that $\gamma\tilde\nu k|B|\deg_H(B)\le
e_G(B,\tilde B\setminus B)$.
We conclude that 
\begin{equation}\label{maxundmoritz}
\maxdeg(H)\le
\frac{\Omega}{\gamma\tilde\nu }=m_\PARAMETERPASSING{L}{lem:sparseRL}\;.
\end{equation} 

\paragraph{Regularising the dense spots in $\tilde B$.} 
We apply Lemma~\ref{lem:sparseRL} with parameters $m_\PARAMETERPASSING{L}{lem:sparseRL},z_\PARAMETERPASSING{L}{lem:sparseRL}$ and $\epsilon_\PARAMETERPASSING{L}{lem:sparseRL}$ as defined by~\eqref{eq:in213} to the graphs
$H_\PARAMETERPASSING{L}{lem:sparseRL}:=\GD$ and
$F_\PARAMETERPASSING{L}{lem:sparseRL}:=H$, together with the ensemble
$\mathcal{B}'$ in the role of the sets $W_i$, and partition of $V(\GD)$ induced by
 $$\mathcal Z_\PARAMETERPASSING{L}{lem:sparseRL}:=\mathcal V\boxplus \big\{V(\Gexp), V(G)\sm V(\Gexp)\big\}\boxplus\big\{V_{\leadsto \smallatoms},V(G)\setminus V_{\leadsto \smallatoms}\big\}\;,$$
 where  $V_{\leadsto \smallatoms}:=\{v\in V(G)\::\:\deg(v,\smallatoms)>b\}$.

Observe that $\mathcal B'$
is an $(\tilde\nu k)$-ensemble satisfying condition~\eqref{lem:sparseRL(item)samesize} of Lemma~\ref{lem:sparseRL},
by~\eqref{eq:ClustersOfRightSize}, by the choice of $k_0$, and by~\eqref{maxundmoritz}.
Thus we  obtain integers $\{p_A\}_{A\in \mathcal{B}'}$ and a family
$\clusters=\{W^{(1)}_A,\ldots,W^{(p_A)}_A\}_{A\in\mathcal B'}$ and a set
$W_0:=\bigcup_{A\in\mathcal B'} W_A^{(0)}$ such that, in particular, we have the
following.
\begin{enumerate}[(I)]
\item We have $\epsilon^{-1}\le p_A\le q_\mathrm{MAXCL}$ for
all $A\in \mathcal{B}'$.\label{aaaaa}
\item We have $|W_{A}^{(x)}|=|W_{B}^{(y)}|$
for any $A,B\in \mathcal{B}'$ and for any $x\in [p_{A}]$, $y\in [p_B]$.\label{bbbbb}
\item \label{ccH}
For any $A\in\mathcal B'$ and any $a\in [p_A]$, there is a set $V\in\V$ for which $W_A^{(a)}\subset V$. We either have that $W_A^{(a)}\subset V(\Gexp)$, or $W_A^{(a)}\cap V(\Gexp)=\emptyset$ and $W_A^{(a)}\subset V_{\leadsto \smallatoms}$, or $W_A^{(a)}\cap V_{\leadsto \smallatoms}=\emptyset$.
\item \label{eq:boundTotalIrregularity}
$\sum_{e\in E(H)}|\mathrm{irreg}(e)|\le \epsilon\sum_{AB\in
E(H)}|A||B|$, where $\mathrm{irreg}(AB)$ is the set of all
edges of the graph $G$ contained in an $\epsilon$-irregular
pair $(W^{(x)}_A,W^{(y)}_B)$, with $x\in[p_A]$, $y\in[p_B]$, $AB\in E(H)$.
\end{enumerate}

Let $\Gblack$ be obtained from $\GD$ by erasing all vertices in $W_0$, and all edges that lie in pairs $(W^{(x)}_A,W^{(y)}_B)$ which are
irregular or of density at most $\gamma^2$. Then
Properties~\ref{defBC:clusters},~\ref{defBC:RL},~\ref{defBC:dveapul}
and~\ref{defBC:prepartition} of Definition~\ref{bclassdef} are satisfied.
Further, Lemma~\ref{notmuchlost} implies~\eqref{eq:sEr}. Together with~\eqref{eq:almostalldashed} we obtain that the number of edges that are not captured
by $( \clusters,\DenseSpots, \Gblack, \Gexp,\smallatoms )$ is at most $(\frac{4\epsilon}{\gamma}+\epsilon\Omega+\gamma+\rho)kn$.

Note that Properties~\eqref{aaaaa}, ~\eqref{bbbbb} and~\eqref{eq:ClustersOfRightSize}
imply that for all $A\in
\mathcal{B}'$ and for any $a\in [p_{A}]$ we have that
\[
\epsilon k\ge |A|\geq |W_{A}^{(a)}|\geq \frac{\tilde\nu k}{q_\mathrm{MAXCL}} \geq \frac{\tilde\nu k}{q}  =\nu k.
\]
Thus also Property~\ref{Csize} of Definition~\ref{bclassdef} holds.

The refinement in~\eqref{eq:Xref} guarantees that the bounded decomposition we have constructed respects the avoiding threshold $b$.

So, it only remains to see Property~\ref{defBC:avoiding} of Definition~\ref{bclassdef}.

\paragraph{The avoiding property of $\smallatoms$.}  
 In order to see
Property~\ref{defBC:avoiding} of Definition~\ref{bclassdef}, we have to show that $\smallatoms$ is $(\Lambda,\epsilon,\gamma,k)$-avoiding with respect to $\DenseSpots$. For this, let $\bar U\subset V(G)$ be such that  $|\bar U|\le \Lambda k$. Let $X$ be the set of those 
vertices $v\in\smallatoms$ that are not contained in any
dense spot $D\in\mathcal D$ for which $|\bar U\cap
V(D)|\le\gamma^2k$. Our aim is to see that $|X|\le \epsilon k$.

Let
$\mathcal{D}_X\subset\DenseSpots$ be the set of all dense
spots $D$ with $X\cap V(D)\neq\emptyset$. 
Setting $\mathcal{A}:=\{A\in \mathcal{\tilde C} :A\cap X\neq\emptyset \}$,
the definition of $\smallatoms$ trivially implies that
$\frac{|X|}{2\tilde\nu k}\le|\mathcal{A}|$. Now, by
the definition of $\mathcal B$, we know that there are at most
$3^{|\mathcal{D}_X|}$ sets $A\in\mathcal A$.
Indeed, for each $D=(U,W;F)\in\mathcal{D}_X$,
 either  $A$ is a subset of $U$, or of $W$, or of
$V(G)\setminus V(D)$. Thus,
\begin{equation}\label{eq_moleculesdetermineatom}
3^{|\mathcal{D}_X|}\geq |\mathcal A|\geq 
\frac{|X|}{\tilde\nu k}\;.
\end{equation}

By Fact~\ref{fact:boundedlymanyspots}, each vertex of $V(G)$ lies in at most $\Omega/\gamma$
of the $(\gamma k,\gamma)$-dense spots from $\DenseSpots$.
Hence
 $$\frac{\Omega}{\gamma}|\bar U|\ge\sum_{D\in\mathcal{D}_X}|V(D)\cap\bar U|\ge |\mathcal{D}_X|\gamma^2 k\overset{\eqref{eq_moleculesdetermineatom}}{\ge}\log_3\left(\frac{|X |}{\tilde\nu k}\right)\gamma^2 k\;,$$ where the second inequality holds by the definition of~$X$. Thus
$$|X|\leq 3^{\frac{\Omega\Lambda}{\gamma^3}}\cdot \tilde\nu k=\epsilon k\;,$$
 as desired.
 This finishes the proof of
Lemma~\ref{lem:decompositionIntoBlackandExpanding}.

\subsection{Sparse decomposition of dense graphs}\label{ssec:sparsedecompofdensegraphs}
Let us explain our remark above that in the setting of a dense graph $G$, Lemmas~\ref{lem:LKSsparseClass} and~\ref{lem:genericBD} produce a regularity partition in the usual sense. So, suppose that $G$ is an $n$-vertex graph and has at least $a n^2$ edges. This needs to be understood with the usual quantification ``$a>0$ is fixed and $n$ is large''.

Recall that when we inquire a $(k,\Omega^{**},\Omega^*,\Lambda,\gamma,\epsilon,\nu,\rho)$-sparse decomposition, the parameters satisfy $\Omega^{**},\Omega^*,\Lambda\gg 1\gg \gamma, \epsilon,\nu,\rho>0$. The interplay between the parameters is quite complicated, and we do not give it here in full (see~\cite[p.\ \pageref{p3.pageref:PAR}]{cite:LKS-cut3} for details). We justify with ``parameter choice'' any further relation we assume between them. Also, let us note that while our exact choice of parameters made in~\cite{cite:LKS-cut3} are tailored for proving Theorem~\ref{thm:main}, we expect these relations to be  satisfied in any application of Lemma~\ref{lem:genericBD}, at least on the loose level we make use of them in this section.

First, we argue that it makes sense to set $k$ linear in $n$, i.e., $k=c n$ for some $c$ depending on $a$ only. Indeed, having $k\gg n$ would allow that all edges of $G$ are uncaptured in~\eqref{eq:generalDecUncapt}, which would make the lemma worthless. On the other hand, with $k\ll n$ we would have all vertices from $Q=\{v\in V(G):\deg(v)>\sqrt{a}n\}$ ending up in the huge-degree set $\HugeVertices$ for which the sparse decomposition provides no structural information. Since $|Q|\ge\sqrt{a}n$, that would be a big loss of information, and thus often undesirable.

So, suppose now that $k=c n$, and suppose that $(\HugeVertices, \clusters,\DenseSpots,\Gblack,\Gexp,\smallatoms)$ is a $(k,\Omega^{**},\Omega^*,\Lambda,\gamma,\epsilon,\nu,\rho)$-sparse decomposition of the dense graph~$G$. Since $\Omega^{**}k=\Omega^{**}cn>n$, we have that $\HugeVertices=\emptyset$. Next, we argue that $\Gexp$ contains no vertices. Suppose on the contrary that it does. Then the minimum degree condition in Property~\ref{defBC:nowheredense}  of Definition~\ref{bclassdef} tells us that $\Gexp$ has at least $\rho k$ vertices of degrees at least $\rho k$ each. Thus, $e(\Gexp)\ge\rho^2k/2=c^2\rho^2n^2/2$. Since $\Gexp$ has at most $n$ vertices, and since $c\rho\gg \gamma$ (parameter choice), we get that $\Gexp$ contains at least one $(\gamma k,\gamma)$-dense spot, a contradiction to $\Gexp$ being nowhere-dense. Last, we claim, that $|\smallatoms|\le \epsilon k$. To this end, consider the set $U_\PARAMETERPASSING{D}{def:avoiding}=V(G)$. We have $|U_\PARAMETERPASSING{D}{def:avoiding}|=n\le \Lambda cn$, and thus the condition in Definition~\ref{def:avoiding} applies. But there cannot exist any $(\gamma k,\gamma)$-dense spot as asserted in Definition~\ref{def:avoiding} since for such a dense spot $D$ we would have $|V(D)|<\gamma^2k$, contradicting its required minimum degree condition. Thus, we conclude that all the vertices $v\in\smallatoms$ are exceptional in the sense of Definition~\ref{def:avoiding}, leading to the desired bound on $|\smallatoms|$.

To summarize, in the sparse decomposition $(\HugeVertices, \clusters,\DenseSpots,\Gblack,\Gexp,\smallatoms)$, we have that $\HugeVertices$, $\Gexp$, $\smallatoms$ are empty or almost empty. Thus, according to Definition~\ref{capturededgesdef}, all the captured edges lie in the regularized graph $\Gblack$.  Property~\ref{bcdef:clustersize} of Definition~\ref{bclassdef} tells us that the clusters have size at least $\nu k=(\nu c)n$, that is, linear in the order of $G$. Further, this property tells us that these clusters are of the same size. We conclude that $\Gblack$ is a regularization of $G$ in the sense of the original regularity lemma.

\subsection{Algorithmic aspects of the decomposition lemma} \label{sssec:DecomposeAlgorithmically}
Let us look back at the proof of the decomposition lemma (Lemma~\ref{lem:decompositionIntoBlackandExpanding}) and observe that we can get a
bounded decomposition of any bounded-degree graph
algorithmically in quasipolynomial time (in the order of the graph). Note that
this in turn provides efficiently a sparse decomposition of any graph, since the
initial step of splitting the graph into huge versus bounded degree vertices 
(cf.~Lemma~\ref{prop:gap}) can be done in polynomial time.

There are only two steps in the proof of
Lemma~\ref{lem:decompositionIntoBlackandExpanding} which need to
be done algorithmically: the extraction of dense spots, and the simultaneous regularization of some dense pairs. 

It will be more convenient to work with a relaxation of the notion of dense spots. We call a graph $H$
\emph{$(d,\ell)$-thick}\index{general}{thick graph} if $v(H)\ge \ell$, and
$e(H)\ge d v(H)^2$. The notion of thick graphs is a relaxation of dense spots, where the
minimum degree condition is replaced by imposing a lower bound on the order, and
the bipartiteness requirement is dropped. It can be verified that in our
proof it is not important  that the dense spots $\DenseSpots$ and the nowhere-dense graph
$\Gexp$ are parametrized by the same constants, i.e., the entire proof would go through even if the spots in $\DenseSpots$ were $(\gamma k,\gamma)$-dense, and $\Gexp$ were
$(\beta k,\beta)$-nowhere-dense for some $\beta\gg \gamma$. Each $(\beta
k,\beta)$-thick graph gives (algorithmically) a $(\beta k/4,\beta/4)$-dense
spot, and thus it is enough to extract thick graphs.

For the extraction of thick graphs we would need to efficiently answer the
following: Given a number $\beta>0$, find a number $\gamma>0$ such that for 
an input number $h$ and an $N$-vertex graph we can localize in $G$ a $(\gamma,h)$-thick graph if it contains a $(\beta,h)$-thick graph, and output NO otherwise.\footnote{We could additionally assume that $\maxdeg(G)\le O(h)$ due to the previous step of removing the set $\HugeVertices$ of huge degree vertices.}
Employing techniques from
a deep paper of Arora, Frieze and Kaplan~\cite{ArFrKa02}, one can solve this problem in quasipolynomial time $O(N^{c\cdot\log N})$. This was communicated to us by Maxim Sviridenko. On the negative side, a truly polynomial algorithm seems to be out of reach, as Alon, Arora, Manokaran, Moshkovitz, and Weinstein~\cite{Alonetal:Inapproximability} reduced the problem to the notorious hidden clique problem, whose tractability has been open for twenty years.
\begin{theorem}[\cite{Alonetal:Inapproximability}]\label{thm:hardness}
If there is no polynomial time algorithm for solving the clique problem for a planted clique of size $n^{1/3}$, then for any $\epsilon\in(0,1)$ and $\delta>0$ there is no polynomial time algorithm that distinguishes between a graph $G$ on $N$ vertices containing a clique of size $\kappa=N^{\epsilon}$ and a graph $G'$ on $N$ vertices in which the densest subgraph on $\kappa$ vertices has density at most $\delta$.\footnote{The result as stated in~\cite{Alonetal:Inapproximability} covers only the range $\epsilon\in(\frac13,1)$. However there is a simple reduction by taking many disjoint copies of the general range to the restricted one.}
\end{theorem}
Of course, Theorem~\ref{thm:hardness} leaves some hope for a polynomial time algorithm when $h=N^{o(1)}$ (which corresponds to $k_\PARAMETERPASSING{L}{lem:decompositionIntoBlackandExpanding}=n_\PARAMETERPASSING{L}{lem:decompositionIntoBlackandExpanding}^{o(1)}$).

\medskip

The regularity lemma can be made algorithmic~\cite{Alon94thealgorithmic}.
The algorithm from~\cite{Alon94thealgorithmic} is based on index pumping-up, and
thus applies even to the locally dense setting of Lemma~\ref{lem:sparseRL}.

\medskip

It will turn out that the extraction of dense spots is the only obstruction to a polynomial time algorithm for Theorem~\ref{thm:main}. In~\cite{cite:LKS-cut3}, we sketch a truly polynomial time algorithm which avoids this step. It seems that the method sketched there is generally applicable for problems which employ sparse decompositions.